\documentclass[10pt,twoside]{article}

\usepackage{fancybox,amsfonts} %MIORTEGA
\usepackage{graphicx}

\usepackage{amssymb}
\usepackage{latexsym}

\makeatletter

\@addtoreset{equation}{section} \makeatother
\newtheorem{theorem}{Theorem}[section]
\newtheorem{remark}[theorem]{Remark}
\newtheorem{lemma}[theorem]{Lemma}
\newtheorem{proposition}[theorem]{Proposition}
\newtheorem{corollary}[theorem]{Corollary}
\newtheorem{definition}[theorem]{Definition}
\newtheorem{example}[theorem]{Example}

\newcommand{\m}{M} %{{\cal M}}
\newcommand{\mo}{M_0} %{{\cal M}_0}

\newcommand{\be}{\begin{equation}}
\newcommand{\ee}{\end{equation}}

\newcommand{\R}{{\mathbb R}}
\newcommand{\LL}{{\mathbb L}}
\newcommand{\SSS}{{\mathbb S}}

\newcommand{\V}{M}

\newcommand{\J}{{\cal J}}

\newcommand{\ben}{\begin{enumerate}}
\newcommand{\een}{\end{enumerate}}
\newcommand{\bit}{\begin{itemize}}
\newcommand{\eit}{\end{itemize}}
\newcommand{\edoc}{\end{document}}
\newcommand{\sm}{\smallskip}

\newcommand{\np}{\newpage}

\newcommand{\bdefi}{\begin{definition}}
\newcommand{\btheo}{\begin{theorem}}
\newcommand{\bprop}{\begin{proposition}}
\newcommand{\brema}{\begin{remark}}
\newcommand{\bcoro}{\begin{corollary}}
\newcommand{\blemm}{\begin{lemma}}
\newcommand{\bexam}{\begin{example}}

\newcommand{\edefi}{\end{definition}}
\newcommand{\etheo}{\end{theorem}}
\newcommand{\eprop}{\end{proposition}}
\newcommand{\erema}{\end{remark}}
\newcommand{\ecoro}{\end{corollary}}
\newcommand{\elemm}{\end{lemma}}
\newcommand{\eexam}{\end{example}}

\title{
\vspace{0.5in} {\bf  Causal boundaries  and holography on wave
type spacetimes }}

\author %{\bf J.L. Flores$^{1}$,
{\bf M. S\'anchez\thanks{The sharp comments by J. L Flores and the
check of a part of the paper by I. R\'acz are warmly acknowledged.
Partially supported by Spanish MEC-FEDER Grant MTM2007-60731 and
Regional J. Andaluc\'{\i}a
Grant P06-FQM-01951. } \\
{\it\small Departamento de Geometr\'{\i}a y Topolog\'{\i}a}\\
{\it\small Facultad de Ciencias, Universidad de Granada}\\
{\it\small Avenida Fuentenueva s/n, 18071 Granada, Spain} }

\begin{document}

%%%%%%%%%%%%%%%%%%%%%%%%%%%%%%%%%%%%%%%%%%%%%MIORTEGA
\setlength{\shadowsize}{2pt}

\newcommand{\scaja}[3]{
\begin{minipage}[c][#1]{14.5em}
\begin{center}\shadowbox{
\begin{minipage}[c][#1]{#2}
#3
\end{minipage}
}\end{center}
\end{minipage}
}

\newcommand{\caja}[3]{
\begin{minipage}[c][#1][c]{10em}\begin{center}
\begin{minipage}[c][#1][c]{#2}
#3
\end{minipage}\end{center}
\end{minipage}
} %%%%%%%%%%%%%%%%%%%%%%%%%%%%%%%%%%%%%%%%%MIORTEGA

\parindent=5mm
\date{}
\maketitle

\begin{quote}

\noindent {\small \bf Abstract.} {\small The notion of {\em causal
boundary} for a spacetime has been a controversial topic during
the last three decades. Moreover, recently the role of the
boundary in the AdS/CFT correspondence for plane waves, have
stimulated its redefinition with some possible alternatives.

Our aim is threefold. First, to review the different classical
approaches to boundaries of spacetimes, emphasizing their
drawbacks and the progressive redefinitions of the c-boundary.
Second, to explain how plane waves and AdS/CFT correspondence come
into play, stressing the  role of the c-boundary as the
holographic boundary in this correspondence. And, third, to
discuss the present-day status of the  c-boundary, making clear
the arguments of a definitive proposal.
}\\

\end{quote}
\begin{quote}

\end{quote}
\begin{quote}

\end{quote}

\newpage
{\small \tableofcontents } \np

%-------------------------------------------------

\section{Introduction}

In the framework of Mathematical Relativity, the  causal boundary,
or c-boundary, is an appealing construction proposed initially by
Geroch, Kronheimer and Penrose \cite{GKP}, in order to attach a
conformally invariant boundary $\partial \V$ to any reasonable
spacetime $\V (\equiv (M,g))$. Roughly, its purpose is  to attach
a boundary endpoint $P$ to any inextensible future or past
directed timelike curve $\gamma$, and the original basic idea is
simple: the boundary point would be represented by the past or the
future of the curve, $P= I^\pm(\gamma)$. By construction, the
c-boundary $\partial M$ is conformally invariant, and does not
contain any direct information on singularities, which may depend
strongly on the conformal factor --$\partial M$ aims to model
``points at infinity'' or even ``conformally invariant
singularities''. Nevertheless, the purely conformal information
contained in $\partial M$ may yield a pleasing picture of the
spacetime. Say, in principle, $\partial M$ would be the union of
three
 disjoint subsets: the {\em future infinity} $\partial_+ M$,
 reached by
 future-directed timelike curves but no past-directed ones, the
{\em  past infinity } $\partial_- M$, dual to the former, and the
{\em timelike boundary} $\partial_0 M$, whose points are reached
by both, future and past directed curves. The latter may represent
naked
 singularities,  the boundary of a removed region in a bigger
 spacetime or, in general,
 losses of global hyperbolicity, i.e., points where
 $J(p,q):= J^+(p)\cap J^-(q)$ is not compact\footnote{Global
 hyperbolicity will be violated  because of these points, whenever $M$
 is  causal, \cite{BeSa2}.}.
 Eventually, asymptotically conformally  flat regions would be represented
 by ``lightlike parts'' in $\partial_\pm M$. Further properties may yield information on black
 holes  or other conformally
 invariant properties. In an ideal scenario, initial conditions for evolution equations
 would be a sort of limit in $\partial_-M$, while boundary conditions would be
 posed on $\partial_0M$.

However, a satisfactory notion of c-boundary must comprise not
only the definition of $\partial M$  but also  the extension of
both, the topology and the chronological relation to the
completion $\bar M = M\cup \partial M$. So, $\bar M$ must satisfy
strong requirements in order to be regarded as satisfactory, and
many specialists have been puzzled along the last three decades
with them. Concretely, with the three apparently harmless
questions:

\begin{itemize}

\item[Q1] {\em Point set definition}. The construction of the
c-boundary $\partial \V$ provides automatically identifications
between the boundary points attached to different (say)
future-directed curves with the same past. Nevertheless, in order
to obtain a  satisfactory c-boundary, some inextensible timelike
curves with {\em opposed time orientations} must also determine
the same boundary points (Fig. %\ref{f1}
1). And a suitable prescription to solve such {\em identification
problem} is not by any means trivial.

\item[Q2] {\em Chronology}. Being the role of Causality obvious
for $\partial M$, causal relations $\ll, < $ %or, at least, the
%chronological one $\ll$ (which is open in $M$),
must be extended to the boundary. Such  extended relations
$\overline{\ll}, \overline{<}$ would also fulfill natural
requirements, as being transitive and equal to the original ones
$\ll, <$ on $M$. Nevertheless, these simple requirements are not
so easy to fulfill, specially when combined with  others (notice
that, if $Q\in \partial M$ satisfies $x\overline{<} Q
\overline{\ll} y$ for some $x,y\in M$, this may yield a new
--spurious-- relation $x\overline{<} y$ or even $x\overline{\ll}
y$, which did not exist for $<, \ll$ in $M$).

\item[Q3] {\em Topology}. The appropriate  topologies for the
boundary $\partial \V$ and the completed space $\bar\V = \V\cup
\partial\V$, should fulfill some natural requirements --for example,
the inclusion $i:M\hookrightarrow \bar M$ must be a topological
embedding and $i(M)$ a dense subset in $\bar M$. But
 subtler questions appear, because there are many simple examples
(say, open subsets of Lorentz-Minkowski $\LL^n$) where both, the
c-boundary as a point set and its topology, seem obvious. It is
not clear at what extent the c-boundary (a conformally invariant
object) must recover such an ``obvious boundary''. Moreover, there
is no a clear rule about which separation properties must satisfy
the boundary points --among them, and  with points in $M$.

\end{itemize}
Relevant contributions by many authors such as Budic and Sachs
\cite{BS}, R\'acz \cite{Ra, Ra2}, Szabados \cite{Sz1,Sz2}, Kuang
and Liang \cite{KLJMP, KLPRD} or Harris \cite{HaJMP,
 HaCQG, HaNonlin}, illuminated both, the possible
definitions of c-boundary and its limitations. Nevertheless, the
notion of c-boundary was unclear at the beginning of the 21th
century, and some authors claimed the impossibility of a natural
general definition.

The interest in the c-boundary has been renewed in the last years,
because of the advances in the holographic principle and AdS/CFT
correspondence. According to this conjectured correspondence, a
string theory on a background spacetime becomes equivalent to a
field theory (the hologram) on its boundary. The question relevant
to us, is which holographic boundary must be chosen. Typically,
the Penrose {\em conformal boundary} is used, as this is the most
common boundary in Mathematical Relativity. Nevertheless, the
detailed AdS/CFT correspondence for plane wave backgrounds
\cite{BMN,BN} stresses the limitations of this boundary. This led
Marolf and Ross \cite{MR1} to consider the c-boundary instead of
the conformal one. In fact, taking into account the progress of
previous authors, they redefined completely the boundary
\cite{MR2}. Remarkably, they introduced the idea of boundary
points as pairs $(P,F)$ of certain past and future sets. Under
this viewpoint, they emphasize that a chronology $\overline{\ll}$
on the completion suggested by Szabados, turns out  consistent
now, and is the unique natural chronology. They also obtained
interesting applications of the c-boundary to the holography of
plane waves, which were independent of the concrete details of
their approximation \cite{MR1}. Moreover, such results were
reobtained and widely extended for the general framework of
wave-type spacetimes in \cite{FS}. But, in spite of this very big
progress, Marolf and Ross breakthrough could not be regarded as
definitive yet. The question was obvious: taking into account the
failure of all previous approaches, to what extent their choices
were unique or undisputable?  For example, in the definition of
the pairs $(P,F)\in\partial M$, Marolf and Ross used Szabados
relation $\sim_S$, which modified a previous one by Budic and
Sachs, is this choice unavoidable? Notice that they introduced two
alternative topologies, and suggested the possible existence of an
intermediate topology with better properties. %Of course, one
%could wonder even if a very different choice would be better.

Recently, Flores revisited systematically the definition of the
c-boundary \cite{F}. Admitting the Marolf-Ross viewpoint that
boundary points must be regarded as pairs $(P,F)$ of some past and
future sets, he introduced a very general notion of {\em
completion} for a spacetime, as well as a natural topology,
inspired in previous work by Harris. Among all the possible
completions, the {\em minimal} --{\em chronological}-- completions
become a privileged (non-empty) subfamily. Then, a detailed study
of the properties satisfied by any completion, is carried out. In
particular, minimal completions satisfy a very reasonable set of
satisfactory properties, required in principle for any c-boundary
(both, as a point set and topologically). The problem of minimal
completions is that, in general, they are not unique. Even though
this problem is not as bad as it sounds --Flores emphasized the
role of chronological completions, and the properties satisfied by
them might be enough for many purposes--, his study goes further.
In general, Marolf-Ross completion is not minimal, but that
completion: (i) shares all the other good properties of the
minimal ones, and (ii) satisfies that any minimal completion can
be naturally included in it. That is, Marolf-Ross completion is
now univocally singled out as a point-set. Moreover, the
properties of Flores topology proved in \cite{F}, plus others to
be discussed here (announcing also results in \cite{FHS}), single
out this topology too. So, Marolf-Ross point-set completion,
endowed with Flores topology,
 and the natural chronology $\overline{\ll}$, are
selected. This choice will be called the {\em Flores boundary}
here, in order to distinguish it from other possible alternatives.
However, the main aim of the present paper is to emphasize that
this must be regarded as the genuine c-boundary, constructed after
the works of all previous authors.

 The present article is organized as follows.
Some brief comments on other boundaries, different to the causal
one, are given in Section \ref{s0}, and their own drawbacks are
also pointed out. The old different approaches for the notion of
c-boundary are reviewed in Section \ref{s1}, especially, questions
Q1, Q2, Q3, are stressed. The roles of holography and plane waves
are emphasized in Section \ref{s2}. In Section \ref{s3} both,
Marolf and Ross \cite{MR1,MR2} and Flores \cite{F} approaches are
explained, and reasons  for our definitive proposal are discussed.
Some of the arguments will be developed further in \cite{FHSst,
FHS}. In Section \ref{s4} we sketch briefly the results and
involved techniques for the boundary of the wave-type spacetimes,
following \cite{FS}. This leads to a highly non-trivial problem
because, in general, the explicit computation of causal boundaries
(TIPs, TIFs) requires new techniques: Busemann functions,
variational interpretations or specific tools for some type of
concrete spacetimes. We point out that the c-boundary of wave-type
spacetimes requires a combination of Busemann-type functions,
variational methods and Sturm-Liouville theory.
% we refer to this reference for the detailed approach. %\footnote{ in this proceedings...
%Que' referencia(s) te parece mejor?}
%for further background on the
%computational problem.
Finally, in Section \ref{s7},  we include a brief summary of our
proposal of c-boundary, for the convenience of the reader.

\section{Preliminaries: some type of boundaries}\label{s0}

In the following, we will consider standard notation and
conventions as in the classical references \cite{BEE, HE, O, W} or
in the recent reviews \cite{GpScqg05, MS}. In particular, any
Lorentzian manifold $(M,g)$ will have signature $(-,+,\dots,+)$,
any spacetime $\V$ is a (connected) time-oriented Lorentzian
$m$-manifold, where the time orientation is assumed implicitly,
and causal vectors are distributed in two cones, each one
containing future or past-directed timelike  ($g(v,v)<0$, $v\in
TM$) and lightlike ($g(v,v)=0$, $v\neq 0$) tangent vectors.

First of all we review some types of boundary used in General
Relativity, different to the causal one. All of them appear under
a natural viewpoint. Nevertheless, we stress that none of them is
fully satisfactory. This will stimulate the efforts to overcome
the
difficulties for %these boundaries and
the c-boundary.

\subsubsection{Penrose conformal boundary.}

In Riemannian Geometry, the stereographic projection yields a
natural conformal open embedding $i:\R^n\hookrightarrow \SSS^n$ of
Euclidean space $\R^n$ into the sphere $\SSS^n$, being the
boundary of the image $\partial i(\R^n)\subset \SSS^n$ a point
which may be interpreted as the conformal infinity of $\R^n$.  In
Lorentzian Geometry there exists also a natural conformal open
embedding $i_P: \LL^n \hookrightarrow \LL^1 \times \SSS^{n-1}$ of
Lorentz-Minkowski $\LL^n$ in the Einstein static Universe
$\LL^1\times \SSS^{n-1}$. As pointed out by Penrose, here the
boundary $\partial i_P(\LL^n)\subset \LL^1\times \SSS^{n-1}$ is
compact and contains certain elements with natural
interpretations: the point $i^0$ or spacelike infinity, the points
$i^\pm$ or timelike infinities and the null hypersurfaces ${\cal
J}^\pm$ or null infinity (see for example, \cite{HE, W}). This
embedding suggests the definition of {\em asymptotically flat
spacetime (at null or spacelike infinity)}, as a spacetime which
admits an open conformal embedding in a bigger ``aphysical''
spacetime $\tilde\V$ (necessarily of the same dimension), which
qualitatively behaves as $i_P$ close to the corresponding elements
at infinity ${\cal J}^\pm, i^0$.

This notion has been widely used in General Relativity: isolated
bodies, mass, or black holes are naturally modelled in
asymptotically flat spacetimes. Nevertheless, the technical
conditions which define an asymptotically flat embedding are very
specific (see \cite[Sect. 11.1]{W} for a discussion). These
conditions imply the essential uniqueness of the asymptotic part
of the boundary (so, this part of the boundary can be regarded as
intrinsic,  \cite{AH}). But the limitations of the approach become
obvious. For a  general spacetime  $\V$ (under, say, some global
reasonable assumption, as being strongly causal or even stably
causal), one can try to find a conformal embedding $i:
M\hookrightarrow \tilde M$ such that (i) $i(M)$ is an open subset
of the ``aphysical'' spacetime $\tilde\V$ (i.e. dim $M$ =dim
$\tilde M$) and (ii) the closure of $i(M)$ in $\tilde\V$ is
compact. Then, the boundary $\partial i(M) \subset \tilde M$ can
be regarded as a sort of conformal boundary. But, in general, this
is neither intrinsic (depends on the embedding) nor systematic
(there is no a way to determine if such conformal embedding exists
and, in this case, how to construct one). In fact, one can check
that a part of $\partial i(M)$ can be regarded as the causal
boundary (even though perhaps with some points artificially
identified). Because of this reason only this part becomes truly
intrinsic and systematic, but even in this case the topology is
not canonical --this will be developed in forthcoming \cite{FHS}.

Further developments in related directions have been carried out.
Garc\'{\i}a-Parrado and Senovilla  introduced {\em isocausal
extensions}, which also yield a conformally invariant boundary
\cite{GpScqg03} --see also \cite{GP-Sa}. Such extensions are less
rigid (and easier to find) than conformal extensions, even though
by this reason their lack of uniqueness is also bigger. %(see \cite{FHS}).
Scott and Szekeres introduced the {\em abstract boundary}
\cite{SS}. Even though their main aim was to study singularities,
this boundary is general, and is defined by using open embeddings
({\em envelopments})  of the underlying {\em differentiable
manifold}. The so obtained a-boundary is unique for the manifold
--as all the possible envelopments are considered. When one
focuses in concrete classes of curves (for example, affinely
parametrized geodesics for a semi-Riemannian metric), some other
boundaries are redefined under this different framework. However,
neither the topology of the a-boundary nor the framework of the
causal boundary were studied in the original article. Even though
some development on the topology was obtained in  \cite{FaSc}, a
further study would be interesting.

\subsubsection{Geodesic and bundle boundaries}

Geroch's g-boundary \cite{GeJMP}  and Schmidt's b-boundary
\cite{Sc} are constructed in order to deal with singularities, so,
they may look  very different to conformal or c-boundaries. The
g-boundary is defined in terms of classes of incomplete geodesics,
and Geroch explored several possibilities for these classes ---the
weakest identifications yield a $T_0$ topology for the quotient.
The b-boundary is a mathematically elegant construction, obtained
by defining a certain positive definite metric on the bundle of
linear frames $LM$ of any semi-Riemannian manifold $(M,g)$; then,
the Cauchy completion of $LM$ induces the b-boundary for $(M,g)$.

Both constructions are systematic, non-conformally invariant, and
satisfy the following {\em a priori desirable} properties :

\bit \item [(i)] every incomplete geodesic in the original
spacetime terminates at a point\footnote{According to Scott and
Szekeres \cite{SS}, a compact manifold should not have neither
singularities nor a boundary. So, as compact Lorentzian manifolds
may be incomplete, this criterium  would not be fulfilled by any
boundary satisfying this a priori desirable condition. },
\item[(ii)] they are geodesically continuous, in a well defined
sense. \eit

The b-boundary was known to give unphysical results for common
solutions to Einstein's equations. For example, the extended
Schwarzschild spacetime $M$ contains a b-boundary point $i$ such
that every neighborhood of $i$ contains all $M$ --and no clear
alternative to Schmidt's construction seemed to exist \cite{Sc2}.

Moreover, a remarkable drawback was found for any construction
satisfying the desirable conditions above: a simple spacetime
$(M,g)$ constructed by Geroch, Liang and Wald \cite{GLW} shows
that minimal conditions (i), (ii) yield a rather undesirable
topology. The example is relatively easy to construct (the
spacetime is just $M=(\R^2\backslash \{s\}) \times \R^2$, $s\in
\R^2$, and the metric $g=\Omega \langle\cdot, \cdot\rangle_{\LL^2}
+ \langle\cdot, \cdot\rangle_{\R^2}$, for some appropriate
function $\Omega>0$; a variant of the example is even flat), and
the a priori undesirable property of the boundary is the existence
of a boundary point $s$ not $T_1$ related to some point in the
spacetime $M$.

Of course, such an example does not affect  the mathematical
validity of the boundaries. So, its true importance, as well as
related constructions of the boundaries, may deserve a further
study (recall, for example, the space of lightlike geodesics
\cite{Lo} and its relation to linking problem \cite{CR}). But, at
any case, such a problem will not affect the causal boundary.

\section{Causal boundaries prior to the holography of waves}\label{s1}

As commented in the Introduction, the purpose of c-boundaries is
to attach a point to any inextensible timelike curve of some
spacetime $M$. As both, the topology and chronological relation
will be extended to the completed spacetime $\bar M$, some mild
causality conditions on $M$ will be used as: {\em chronology}
(inexistence of closed timelike curves), {\em causality}
(inexistence of closed causal curves), to be {\em (past or future)
distinguishing} (future distinguishing: different $p,q\in M$ have
different $I^+(p), I^+(q)$; past distinguishing: analogous with
$I^-(p), I^-(q)$; distinguishing: future and past distinguishing)
or {\em strong causality} (inexistence of ``almost closed'' causal
curves). Recall that strong causality is also equivalent to the
equality between the natural topology of $M$ and its Alexandrov
topology, i.e., the one generated by $I^\pm(p)$ for all $p\in M$.

\subsection{Starting c-boundaries: GKP construction}
\label{s3.1}

The Geroch, Kronheimer and Penrose \cite{GKP} boundary is a
general construction, in principle applicable to any
 strongly causal spacetime $\V$, and it is explicitly
 intrinsic, systematic and unique.

The construction has a very appealing first part: the definition
of the {\em precompletion} $\V ^\sharp$, which, essentially, is
retained in all the subsequent developments of the c-boundary.
Briefly (see for example \cite{Ha09} in this proceedings for more
details), one starts by declaring that $P\subset M$ is a past set
if $P =I^-(P)$ and, then, $P$ is an IP (indecomposable past set)
if it cannot be written as the union of another two past subsets.
Such an IP is necessarily either {\em proper} (PIP) or {\em
terminal} (TIP). In the former case, $P=I^-(p)$ for some $p\in M$
(so, $M$ itself is identifiable to the set of all PIP's, as the
spacetime is past distinguishing), in the latter, $P=I^-(\gamma)$
for some inextendible future-directed timelike curve

The set of all TIP's is the {\em future preboundary}
$\hat{\partial} M$ of $M$, and the set  of all the IP's is the
future precompletion $\hat{M}$. Analogously, indecomposable future
sets, (IF's), which may be either proper (PIF's) or terminal
(TIF's) are considered, and one defines the past preboundary
$\check\partial{M}$ and past precompletion $\check{M}$. The
precompletion $M^\sharp$ of $M$ is essentially  $M$ plus the
preboundary points or, more precisely: $M^\sharp = (\hat{M} \cup
\check{M})/\sim$ where $\sim$ is the relation of equivalence
$I^+(p)\sim I^-(p), \forall p\in M$.

However, the next steps in the GKP construction concerns the
questions Q1, Q3 in the Introduction (identifications and
topology), and they became widely controversial. The GKP
construction tries to solve both questions at the same time. First
$M^\sharp$ is topologized by taking as a subbase the sets
$F^{int}$, $F^{ext}$, $P^{int}$, $P^{ext}$:
%where, say, for the ``$F$'' case:
$$F^{int}= \{P\in \hat M: P\cap F\neq
\emptyset\}, \quad , \quad F^{ext}= \{P\in \hat M: P=I^-[\omega]
\Rightarrow I^+[\omega]\not\subset F\},  \quad \quad \forall F\in
\check{M},$$ and analogously for $P\in \hat M$. Notice that, when
$F$ is a PIF, $F=I^+(p)$, then the set of PIP's corresponding to
$F^{int}$ is equal to $F$ itself, and  the set corresponding to
$F^{ext}$ is $M\backslash \bar J^+(p)$. In this sense, this
topology is  Alexandrov-type.  %(i.e., the topology of $M$ generated
%by the sets type $I^\pm(p)$, which coincides with the natural
%topology of $M$ in strongly causal spacetimes).
As pointed out in \cite[Figure 6]{GKP}, the sets type $P^{ext},
F^{ext}$ are required to give a reasonable basis for the topology
in the possible ``lightlike part'' of the boundary. Finally, the
causal completion would be defined as the quotient $\bar M=
M^\sharp/R_H$, where $R_H$ is the minimum identification of {\em
preboundary points} to obtain a Hausdorff space.

\sm\sm

\noindent This {\em a priori} imposition of Hausdorffness was
found unsatisfactory by several authors. Perhaps, the most
surprising drawback was pointed out by Kiang and Lian
\cite{KLJMP}: the GKP boundary recovers well the natural boundary
of a timelike half plane of $\LL^2$, but this is not the case for
a
half space of $\LL^3$ (Figs. 1 %\ref{f1}
and 2). %\ref{f2}).
However, this was not the unique problem:

(1) One such identification only between preboundary points
yielding a Hausdorff quotient (and, then, the intersection of all
them, $R_H$) may not exist in general, see  \cite{Sz1}.
Intuitively, the reason is that strong causality might be lost
``at a boundary point'' and, so, one such point might be
non-Hausdorff related to points of the spacetime. Nevertheless,
such an identification does exist if the spacetime is stably
causal; so, the approach must be restricted to this class of
spacetimes.

(2) Taub spacetime is static, but its GKP ``singularity'' is only
a point, not a  line \cite{KLLPRD}, as one would expect
--unsatisfactory properties for plane waves under GKP and other
approaches will be also found  \cite[Section 5]{MR2}.

 (3) For $\LL^n$, the topology of the (causal part of the) Penrose
conformal boundary does not agree with the GKP topology
\cite{HaCQG}\footnote{\label{footh} The boundary $\hat\partial
\LL^n$ is a cone, as in Penrose conformal embedding. But each one
of its null generators is a GKP-open subset, by arguments similar
to
those in Fig. 2. %\ref{f2}.
Such a problem will appear also in other approaches which
essentially maintain the GKP topology (see also \cite{FH}).}.

\begin{figure}[ht]\label{f1}
\begin{center}
\includegraphics[width=.5\textwidth]{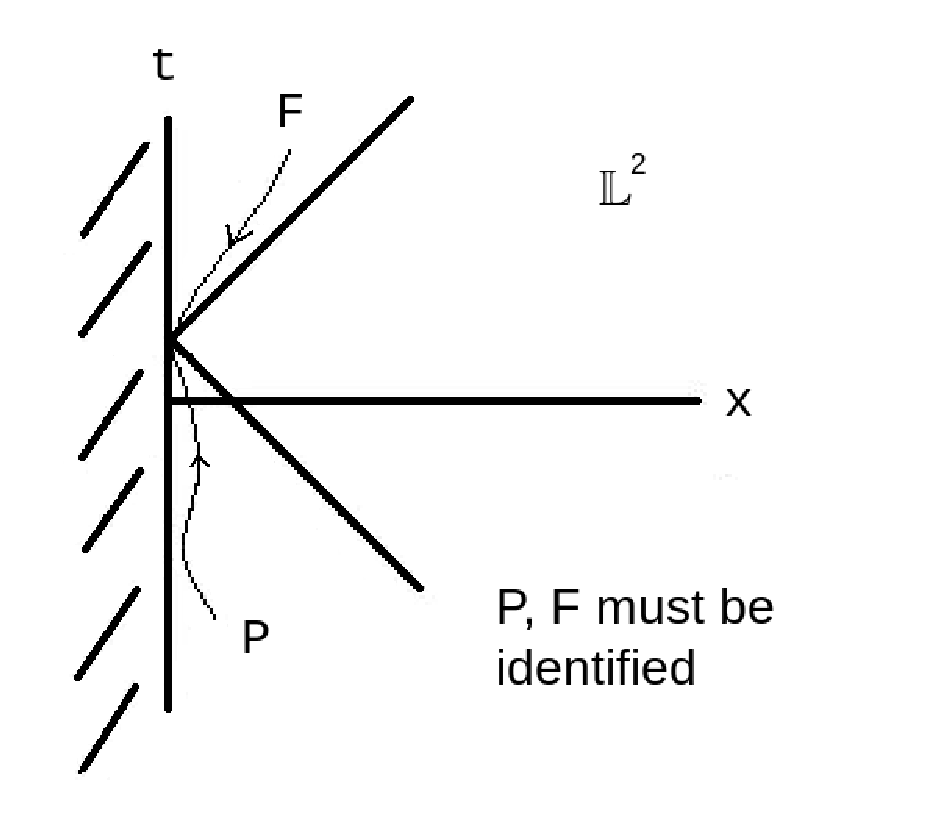}
\caption{Spacetime $\LL^1\times \R^+ \subset \LL^2$ ($x>0$). The
relation of equivalence $R_H$ must be introduced in order  to
recover the natural boundary: as $P=I^-(\rho), F=I^+(\gamma)$
``finish'' at the same point with $x=0$, they have to be
identified to a single point $Q$. This happens in the GKP topology
as $P$, $F$ are $T_1$-related ($P\in F^{ext}$, $F\in P^{ext}$) but
not $T_2$-related in $M^\sharp$ (say, $P^{ext}\cap F^{ext} \neq
\emptyset$).}
%\medskip
%\parbox{5cm}{
%\footnotesize
%{\small Conformal embedding of Minkowski
% spacetime in the Einstein static universe.}

\end{center}

\end{figure}

\begin{figure}[ht]\label{f2}
\begin{center}
\includegraphics[width=0.8\textwidth]{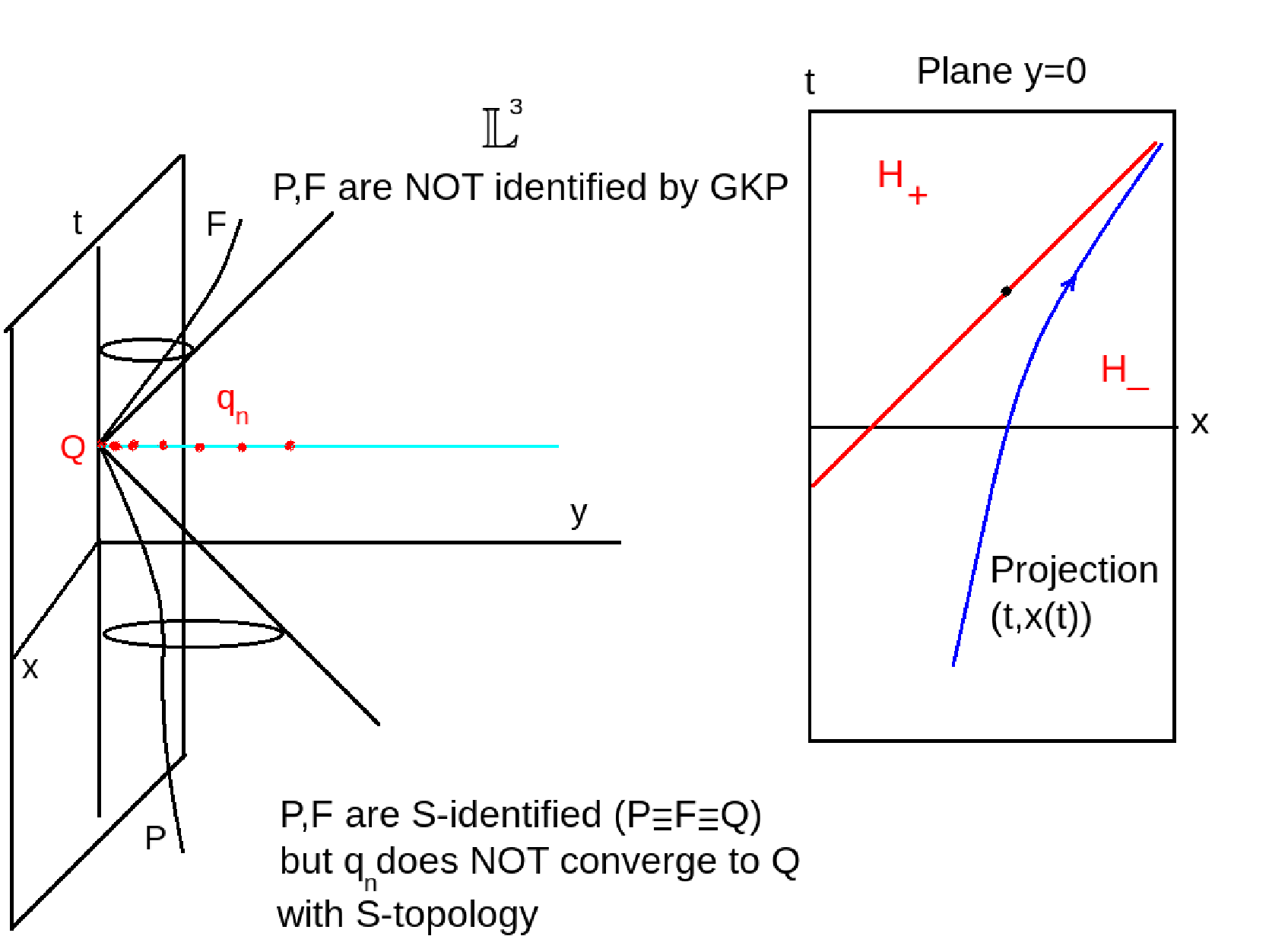}

%\medskip
%\parbox{5cm}{
%\footnotesize
%{\small Conformal embedding of Minkowski
% spacetime in the Einstein static universe.}

\caption{ GKP construction does not recover the natural boundary
of $\LL^2\times\R^+\subset \LL^3$ ($y>0$): $P=I^-(\rho),
F=I^+(\gamma)$  are not identified as they are $T_2$-related. In
fact, a lightlike plane which intersects orthogonally the
(conformal) boundary $y=0$ at $Q$ % corresponding to $P,F$
 yields two hemi-spaces $H_\pm$. Now,
 $H_- =I^-(\alpha)$ for some timelike $\alpha(t)$ (type $\alpha (t)=(t,x(t),(1+t^2)^{-1})$). So, $H_-$  is a TIP, $H_+$ a TIF, and $H_\pm^{ext}$ separates $P$ and
$F$ ($P\in H_+^{ext}$, $F\in H_-^{ext}$ and $H_+^{ext}\cap
H_+^{int}=\emptyset $), \cite{KLJMP}. Szabados relation does
identify $P,F$ to the single point $Q$, but  $\{q_n\}_n
\not\rightarrow Q$ with his topology \cite{KLPRD}.}
\end{center}
\end{figure}

 \subsection{First modifications}

 Budic and Sachs \cite{BS}, R\'acz \cite{Ra}
 and Szabados \cite{Sz1, Sz2}, proposed some modifications of the GKP
construction, in order to overcome previous problems. Let us start
with the second one.

\subsubsection{R\'acz's topology}

The simplest modification is  to change the ``Alexandrov type''
topology. R\'acz proposed the modifications of the subsets
$F^{int}$, $F^{ext}$ (and analogously for $P$). Essentially: (i)
these open sets are defined
 only when $F$ is a PIF (no a TIF), and (ii) they generate not only IPs but also IFs . More precisely:
$$F^{int}= \{A\in \hat M \cup \check M: A\in \hat M \;\hbox{and}\; A\cap F\neq
\emptyset \; \hbox{or} \; A\in \check M \;\hbox{and}\; A=I^+[S]
\Rightarrow I^-[S] \cap F \neq \emptyset \},$$ % \quad , \quad
$$F^{ext}= \{A\in \hat M \cup \check M: A\in \check M
\;\hbox{and}\; A\not\subset F \; \hbox{or} \; A\in \hat M
\;\hbox{and}\; A=I^-[S] \Rightarrow I^+[S]\not\subset F \},$$

Some problems of the GKP construction do not appear in R\'acz one
(for example, the boundary of Taub is 1-dimensional --but again
the spacetime must be  stably causal, in order to ensure the
existence of the relation of equivalence $R_H$). However, Kuang
and Liang \cite{KLPRD} found an example unfavorable to R\'acz
topology (Fig. 3). %\ref{f3}).

\begin{figure}[ht]\label{f3}
\begin{center}
\includegraphics[width=0.5\textwidth]{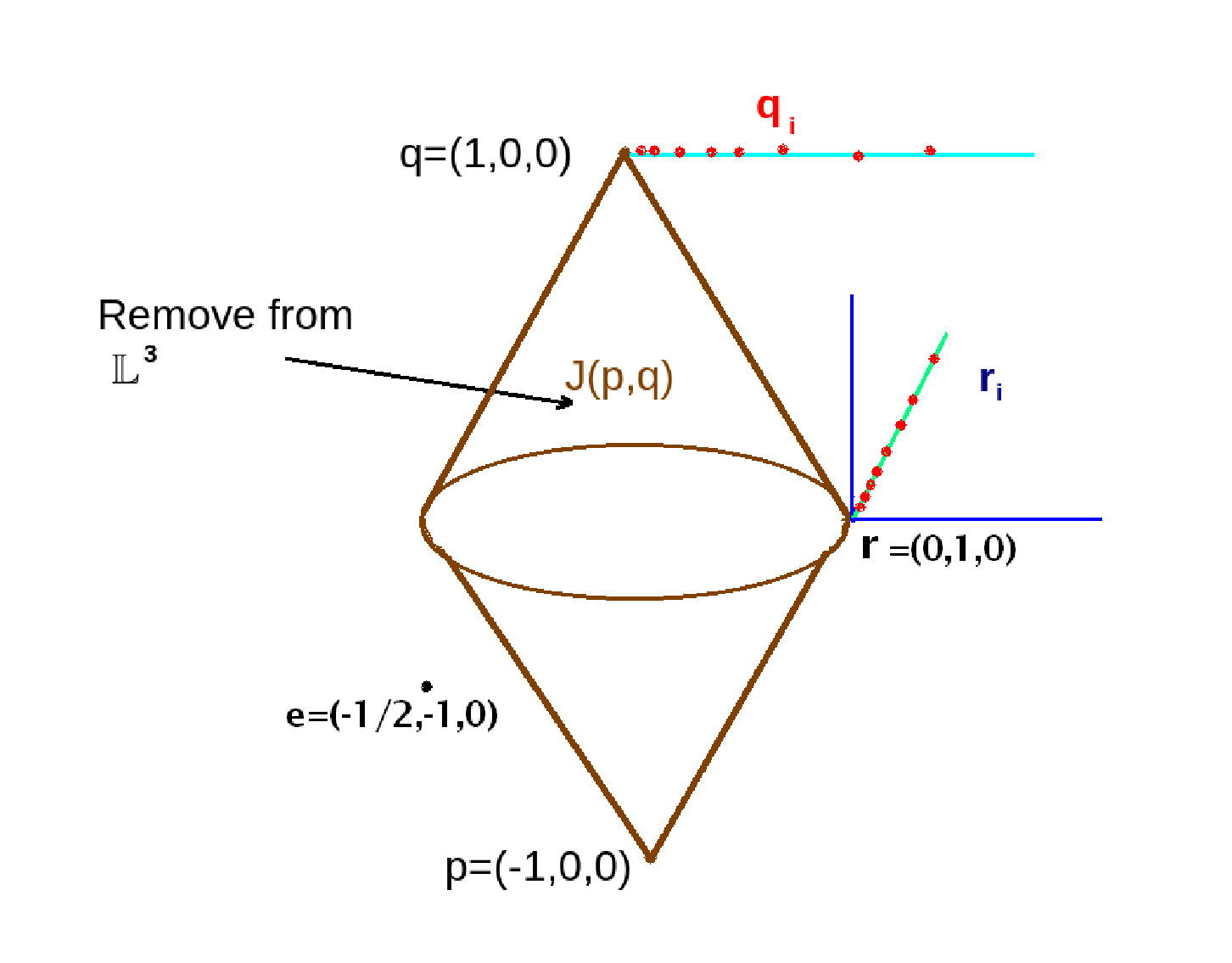}

\caption{Causally continuous spacetime $M=\LL^3\backslash
(J^+(t=-1,x=0,y=0)\cap J^-(1,0,0))$. The removed points
$r=(0,1,0), q=(1,0,0)$ are naturally identified with boundary
points. Nevertheless: (a) for R\'acz topology, $\{r_n=(0,1,1/n)\}
\not\rightarrow r$, as $I^+(-1,2,0)^{int}\cup I^-(1,2,0)^{int}$
projects onto an open subset of the completed space which contains
$r$ but no $r_n$ \cite{KLPRD}, (b) for  BS completion,
$\{q_n=(1,1/n,0)\}\not\rightarrow q$, as $e=(-1/2,-1,0)$ is
 not chronologically related to any $q_n$,
but $e\ll q$   in the completed spacetime $\bar M$ --that is,
$I^+(e)$ yields an open subset of $\bar M$  containing $q$ but no
$q_n$. }
\end{center}
\end{figure}

Notice that R\'acz topology on $M^\sharp$ maintains the  GKP
relation $P\in F^{ext}$ pointed out in Fig. 1. %\ref{f1}.
In principle, this behavior seems rather artificial, but becomes
essential  in order to ensure that $P, F$ are not $T_2$ related
--and, thus must be identified with the same point.

\subsubsection{Budic and Sachs' equivalence relation}

Budic and Sachs  made a full revision of the GKP construction. One
of the main ingredients they introduced (which will turn out to be
essential for our final proposal) is a direct identification
between TIPs and TIFs --recall that both, GKP and R\'acz
identifications were defined indirectly. Namely, for $P\in \hat M,
F\in \check{M}$:
$$ P\sim_{BS}F \quad \Longleftrightarrow \quad P=\downarrow F %(:= \hbox{Int} \{I^-(x): x\in F\})
\;\; \hbox{and} \;\; F=\uparrow P$$ where, say, the common past
$\downarrow F$ is the interior of $\cap_{x\in F}I^-(x)$. Then,
they also  extended both chronological and causal relations
 on the quotient space $\overline{M}=M^\sharp/\sim_{BS}$. By using these relations,
 an Alexandrov-type topology was defined on $\bar M$. The
 construction seemed specially good for {\em causally continuous}
 spacetimes, but  Kuang and Liang \cite{KLJMP} also found a  unfavorable example in this case (Fig.
 3). % \ref{f3}).
%They showed that
% the spacetime above $M=\LL^3\backslash J(p,q)$
%$p=(-1,0.0), q=(1,0,0)$  is , and  the sequence
Intuitively,  the BS approach is conceived for spacetimes $M^m$
with no ``$(m-1)$-dimensional'' parts removed, or
``$m$-dimensional holes". In fact, if one removes a $m-1$
dimensional part,
undesirable properties may appear, see Fig. 4. % \ref{f4}.
Nevertheless, causal continuity cannot prevent the existence of
holes, as Fig. 3 % \ref{f3}
shows. Moreover, analogous problems for causally simple spacetimes
were also pointed out by R\"ube \cite{Ru}.
%\footnote{Ese ejemplo no me parece muy claro. Sin embargo, puede
%ser un buen test para Flores topology}.
The problems of the identifications become then specially
delicate, as causal simplicity is the step in the standard causal
hierarchy of spacetimes which follows to causal continuity, and
the next step is global hyperbolicity
---where the problem  becomes trivial, as no identifications can appear.

\begin{figure} \label{f4}
\includegraphics[width=.5\textwidth]{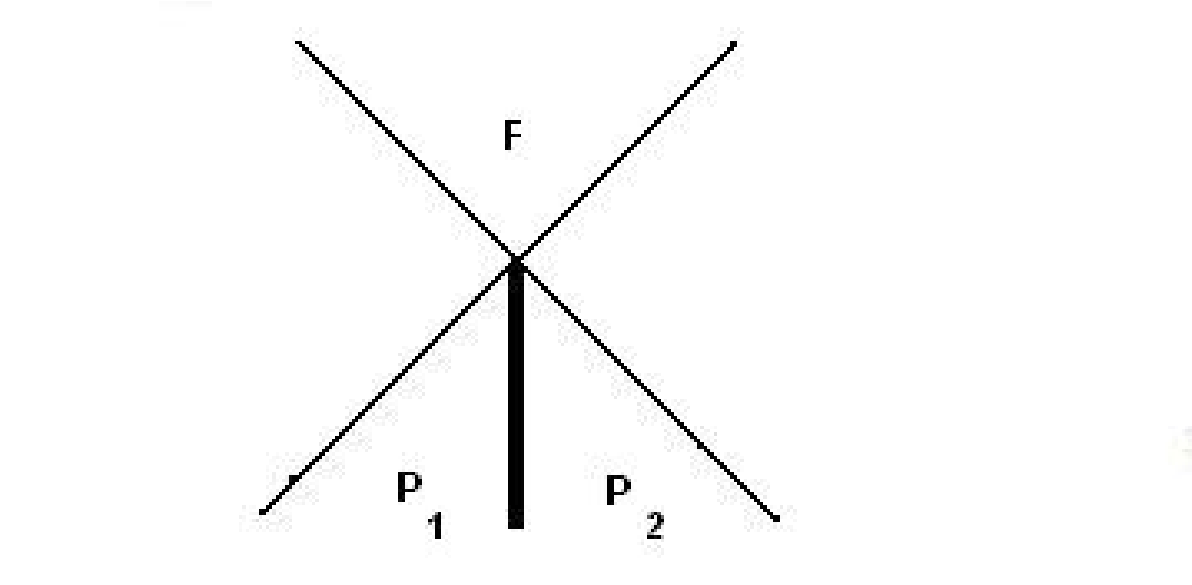}
%\medskip
%\parbox{5cm}{
%\footnotesize
\caption{ Because of the removed segment  of $\LL^2$
%$(-1,1)\times(-1,1)\subset \LL^2$,
$\downarrow F= P_1 \cup P_2$. Thus $P_i \not\sim_{BS} F $ for
$i=1,2$. As $P_1\sim_S F \sim_S P_2$, Szabados imposes that the
three sets represent a single boundary point. Marolf-Ross will
take both pairs $(P_i,F)$ as distinct boundary points.}
\end{figure}

\subsubsection{Szabados reformulation}

Szabados made a penetrating study on which points of the
GKP-boundary can be separated from points of the spacetime. As
stable causality or causal continuity do not prevent the
appearance of undesirable situations, he tried to give a solution
for any strongly causal spacetime.

Szabados reformulated the BS identification in $M^\sharp$, in
order
to deal also with examples as in Fig. 4. %\ref{f4}.
In the set of all the IPs and IFs $\hat M \cup \check M$, he
introduces the relations:

\be \label{eSz} P\sim_S F \Longleftrightarrow \left\{
\begin{array}{l} F \quad
\hbox{is included and  is maximal in} \quad \uparrow P \\
P \quad \hbox{is included and is maximal in} \quad \downarrow F
\end{array} \right. \ee
and extend this by transitivity to a relation of equivalence (here
``maximal'' means with respect to the partial ordering $\subset$).
Moreover, he also introduced a second technical identification.
This involves terminal indecomposable sets of the same type
(TIP-TIP or TIF-TIF), as sometimes this becomes natural (see Fig.
5). %\ref{f5}).

\begin{figure}\label{f5}
\includegraphics[width=.4\textwidth]{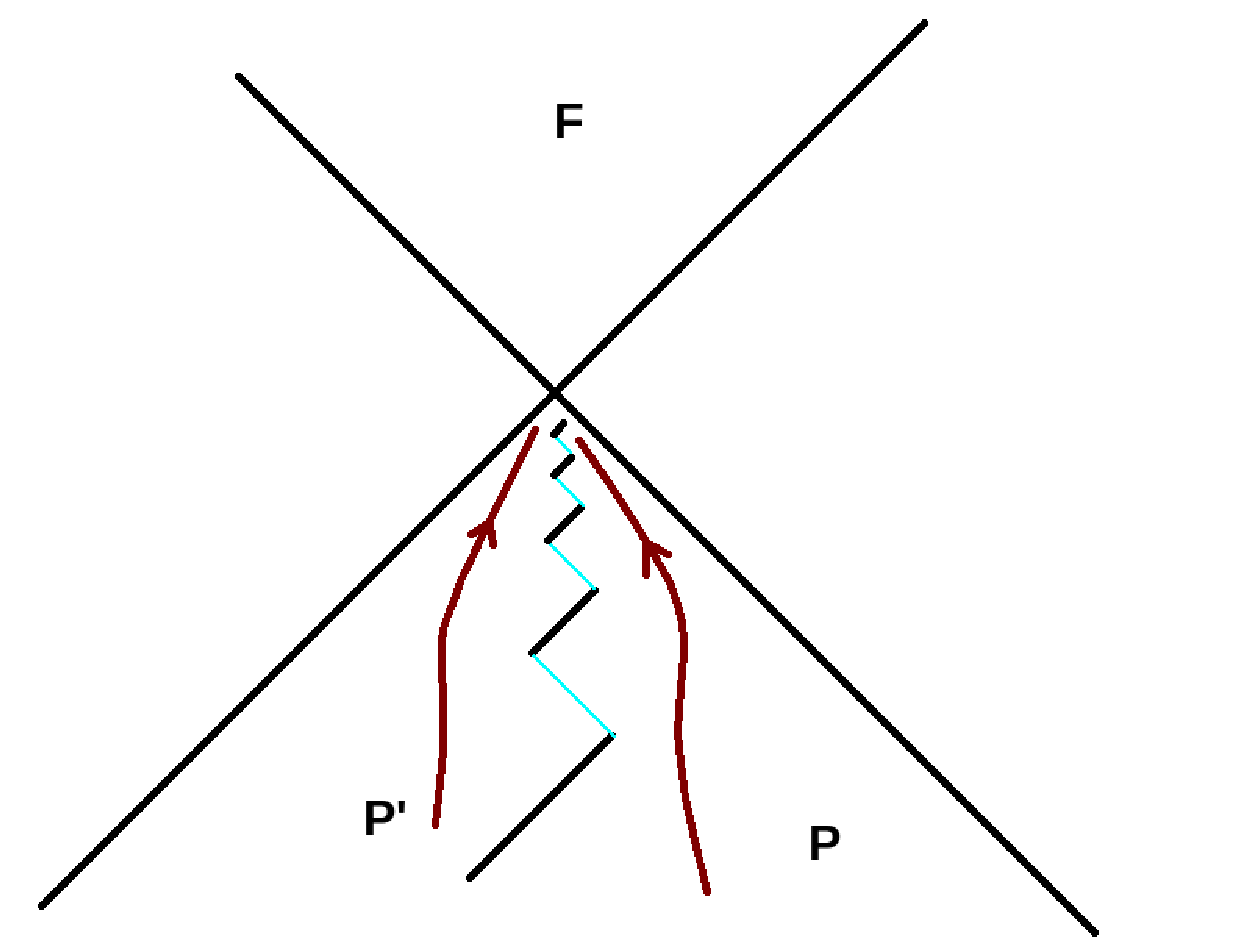}
%\medskip
%\parbox{5cm}{
%\footnotesize
\caption{ In $\LL^2$ minus a removed sequence of thick segments
(at 45 degrees) one has: $P'\subsetneq P$, $F=\uparrow P= \uparrow
P'$; thus $P'\not\sim_S F$. A second Szabados relation identifies
$P$ and $P'$. Marolf-Ross take both pairs $(P,F), (P',\emptyset)$
and their topology is not $T_1$ here.}
\end{figure}

\sm This defines the boundary as a point set.  A chronological
relation  and a  topology were also introduced. Even though the
chronological relation had some problems inherent to the
identification approach pointed out by Marolf and Ross (see
Example \ref{exfl}(3) below), it will suggest the final choice of
the chronological relation. About the topology, some previous
problems of the GKP boundary do not appear, but Kuang and Liang
\cite{KLPRD} found a new very unfavorable example. In fact, the
Szabados construction does not recover well the topology for the
completion of $\LL^2\times\R^+$, Fig. 2. %\ref{f2}.

What is worse, the Kuang and Liang example will hold whenever: (i)
the GKP topology is used and (ii) the identifications rule is such
that a TIF $F=I^+(\lambda)$ and a TIP $P=I^-(\gamma)$
corresponding to the boundary of the removed region $y\leq 0$ are
identified iff $\lambda, \gamma$ have the same point $(t_0,x_0,0)$
as their endpoint. Therefore, these authors claimed \cite{KLPRD}:

\begin{quote}
{\em We are inclined to believe that the whole project of
constructing a singular boundary has to be given up.}
\end{quote}
Even though there were reasons for this pessimistic
viewpoint\footnote{\label{foot}Apart from Kuang and Liang
counterexamples, Harris not only pointed out that the problems of
the topology of $\partial \LL^n$ were inherent to GKP-type
constructions (see footnote \ref{footh}), but also gave an example
which suggested that the identifications may yield ``incomplete
completions'', see \cite[Appendix]{HaJMP}.}, the analysis of
previous approaches makes reasonably clear: (a)  GKP
identification rule must be abandoned, (b) GKP topology (and its
naive variations) have to be changed, and (c) one cannot forget to
find a good extension of
$\ll$ to the boundary. The %The first item had been
Budic-Sachs and Szabados approaches introduced modifications in
the directions (a) and (c). But for deeper changes, including (b),
 it is convenient to study first a
simplified case, free of bothering identification rules.

\subsection{ Harris'
universal construction }

Harris \cite{HaJMP, HaCQG} (see also updated \cite{HaGeoD})
focused on the less problematic future part $M, \hat M, \hat
\partial M$ or {\em future chronological (chr) boundary}. This is useful
from the practical viewpoint: in many cases one may be interested
only in what happens ``towards the future'' (or the past) with no
annoying discussions about identifications --and in spacetimes
such as the globally hyperbolic ones, this always will happen. But
this viewpoint will be also useful to attain a general definition
of c-boundary: on one hand, the universal properties of the
chronological boundary will give a support for the c-boundary, on
the other, the natural topology for the chr boundary, will suggest
a natural topology for the c-boundary.

As Harris focuses on $\hat \partial M$ as a point set, the first
task will be to extend $\ll$ to all $\hat \partial M$. He  does
not try to extend the causal relation $\leq$, and we will not
worry about it  (we will give some brief comments on extended
$\leq$ in Subsection \ref{s5.5}).

\subsubsection{The chronological relation in
$\hat\V$}\label{3.3.1}

Let $M$ be a strongly causal spacetime, and $\hat M$ be its future
GKP (pre)completion.  %In general, $\hat M$ will not admit a smooth
%structure, but extensions of both, the $M$-topology and
%chronological relation $\ll$ would be required. Starting by the
%latter, the first problem is: which  extension  of $\ll$ must be
%chosen?
Following GKP approach,  Harris defined a chronological relation
 extended to the completion $\hat M$, $\ll^c$, as: \be \label{eell}
\begin{array}{lll}
x\ll^c y &\Leftrightarrow & x\ll y,\\

x\ll^c Q &\Leftrightarrow & x \in Q,\\

P\ll^c y &\Leftrightarrow & P \subset I^-(z) \; \hbox{for some} \; z\ll y\\

P\ll^c Q &\Leftrightarrow & P \subset I^-(y) \; \hbox{for some} \;
y\in Q
\end{array}\ee
for any  $x,y\in\V$, $P,Q\in \hat\partial\V$. As Harris points
out, $\ll^c$ becomes nice when $\V$ is {\em past-determined},
i.e.: $x\ll z$ holds if $I^-(x)\subset I^-(y)$ for some $y\ll z$.
Otherwise,  the following difficulty appears. Consider $x\in \V$
and its past $P=I^-(x)$. Because of the different cases in the
definition of $\ll^c$, perhaps $x\not\ll^c y$ but, when one
removes $x$ from the spacetime, $P\ll^c y$ (Fig. 6). %\ref{f6}).

\begin{figure}\label{f6}
\includegraphics[width=1\textwidth]{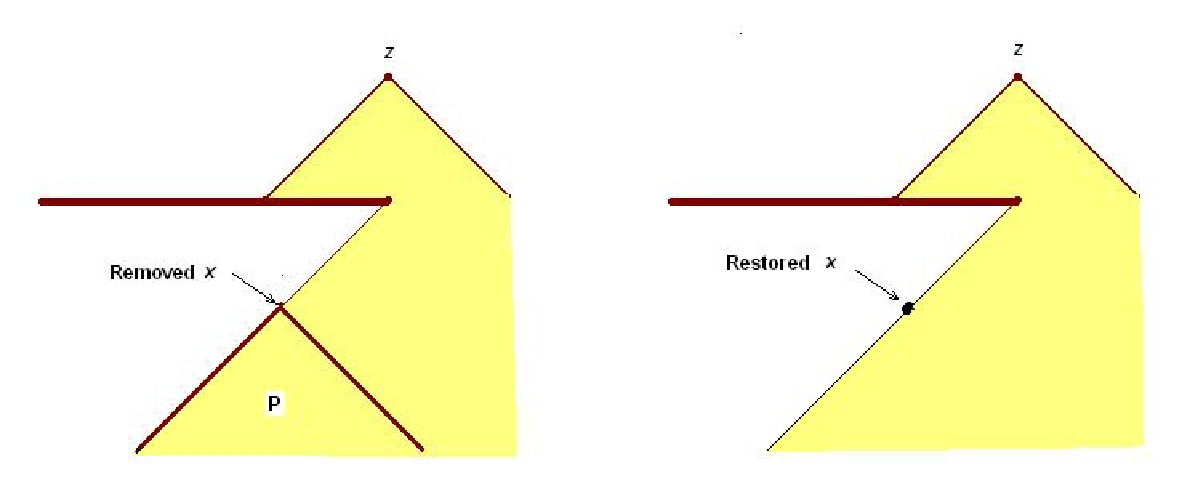}

\caption{On the right, $Y=\LL^2\backslash\{(0,x): x\leq 0\}$; on
the left, $X$ is equal to $Y$ with the point $x$ removed ($X, Y$
are not past-determined). The natural inclusion $i:X\rightarrow Y$
cannot be extended to a chronological map $\hat {\i}: \hat
X\rightarrow \hat Y$ if one considers the  chronological relation
$\ll^c$ in  the completions. In fact $P\ll^c z$ but $x=\hat {\i}
(P) \not\ll^c z$. }
\end{figure}

In order to overcome this difficulty, Harris will also consider a
{\em past-determined chronological relation} $\ll^P$ in $\hat\V$,
which can be reformulated as follows:
\begin{equation}\label{extendedll} P\ll^P Q \Leftrightarrow \uparrow
P \cap Q \neq \emptyset, \quad \quad \hbox{for any IPs}\; P,Q \in
\hat\V . \end{equation} Notice that $\ll^P$ considers both, TIPs
and PIPs on the same footing. Nevertheless, $\ll^P$ introduces new
(spurious) relations in the original spacetime $\V$ --so that it
will be past-determined. Sometimes, this difficulty will make
necessary  to focus on the class of past-determined spacetimes
(compare with Remark \ref{rh} below).

\subsubsection{Category of chronological sets}

Next, one has to construct a category of spaces with includes
both, strongly causal spacetimes and spacetimes-with-boundary.

\bdefi \label{dobj} A {\em chronological set} $(X,\ll)$ is a  set
$X$ endowed with a binary relation $\ll$ which is transitive,
anti-reflexive and satisfies:

(i) it contains no isolates: each $x$ satisfies $x\ll y$ or $y\ll
x$ for some $y\in X$, and

(ii) it is {\em chronologically separable}, that is, there exists
a countable ${\cal S}\subset X$ which is {\em chronologically
dense}: for all $x \ll y $ there exists $s\in {\cal S}$ such that
$ x\ll s \ll y$.\edefi

All these conditions are satisfied trivially in chronological
spacetimes. About the inexistence of isolates, recall that in
spacetimes each $x$ satisfies both $x\ll y$ and $y'\ll x$, for
some $y, y'$. Nevertheless, to impose only one of them, becomes
natural for points of the boundary. The chronologically dense
subset ${\cal S}$ will ensure appropriate technical properties.

Notice that we have not defined a topology in $(X,\ll)$ yet.
However,  definitions for spacetimes on causality and GKP
construction ($I^\pm(x)$, IPs, TIPs,...) can be translated
immediately. Future-directed timelike curves are replaced
naturally by {\em future-directed {\em (or} chronological{\em )}
chains}: sequences $c=\{x_j\}_j$ such that $x_j \ll x_{j+1}$ for
all $j$. Chronological separability ensures that IPs coincide with
pasts of future-directed chains.

\bdefi A point $x\in X$ is a {\em future limit} of the
future-directed chain $c$ if $I^-(c)=I^-(x)$. $X$ is called {\em
future-complete} if any such future-directed chain  has a future
limit in $X$.\edefi Future-directed chains provide a closer
analogy with the Cauchy completion of a metric space, as they play
a role similar to Cauchy sequences.

\bdefi The {\em chronological boundary} $\hat \partial X$ of
$(X,\ll)$ is the set of all the TIPs of $X$, and the {\em
chronological completion} $\hat X$ is the union $X\cup \hat
\partial X$. The latter is also regarded as a chronological
space %\footnote{According to Harris,
with the extended relation $\ll^c$ %which may be chosen as $\ll^e$ or $\ll^c$
in (\ref{eell}). %--recall also Remark \ref{rh}.}.
\edefi
 %At any case, both coincide if $(X,\ll)$ is
%{\em past-determined} in the same sense explained above.
Due to the generality of chronological spaces, there are two
conditions which will be always imposed and are automatically
satisfied in both, strongly causal spacetimes and their
completions: \bit \item {\em Past-distinguishing}: $x\neq y
\Rightarrow I^-(x)\neq I^-(y)$. This allows to identify each point
$x\in X$ with its PIP, $I^-(x)$, and  ensures the uniqueness of
the future limit for chronological sequences. Notice that strong
causality is not considered for $(X,\ll)$: first, a topology is
needed for its definition, but even when this is done, a strongly
causal spacetime might have, in principle,  boundary points where
such strong causality were violated \cite{Sz1}.%\footnote{Comprobar}
\item {\em Past-regularity}: $I^-(x)$ is an IP for all $x\in X$.
This condition is completely natural, but one might be inclined to
drop it in same cases. For example, remove the negative time axis
of $\LL^2$, and consider its chronological completion $\hat X$
(see Fig. 4). %\ref{f4}).
Each removed point $(t,0), t\leq 0$
yields naturally two boundary points $(t,0^-), (t,0^+)$. Then,
identifying $(0,0^-)$ and $ (0,0^+)$, one obtains a natural {\em
non} past-regular chronological set. \eit A natural morphism
between chronological spaces $X, Y$ will preserve the chronology.
But as future limits become fundamental for chronological
completions, we will consider morphisms which also preserve these
limits:

\bdefi \label{dmorf} A function $f: X\rightarrow Y$ between
chronological sets is {\em future-continuous} if it:

(i) is chronological: $x_1\ll x_2 \Rightarrow f(x_1) \ll f(x_2)$,
and

(ii) preserves future limits: if $x$ is the future limit for a
future-directed chain $\{x_j\}_j$ then so is $f(x)$ for the
(necessarily future) chain $\{f(x_j)\}_j$. \edefi

\begin{remark}\label{rhfc} {\em If $X,Y$ are strongly causal
spacetimes (and, thus, endowed with  Alexandrov topology, which
comes directly from chronology) a function is future and past
continuous iff it is continuous and carries timelike curves to
timelike curves \cite[p. 5433s]{HaJMP}. Nevertheless,
future-continuity does not ensure continuity (consider in $\LL^2$
the map $(t,x)\mapsto (t-1,x)$ for $t\leq 0$ and equal to the
identity otherwise).}
\end{remark}

Definitions \ref{dobj}, \ref{dmorf} allow to consider the category
${\cal C}$ with objects the past-distinguishing past-regular
chronological sets and morphism the future-continuous functions.
As the morphisms preserve completeness, we can also consider the
subcategory ${\cal C}_0$ whose objects are all the future-complete
${\cal C}$-objects. The next step in Harris approach is to show
that every morphism $f: X\rightarrow Y$ extends to a morphism
$\hat f: \hat X \rightarrow \hat Y$ in order to arrive at a
functor $ \widehat \; : {\cal C}\rightarrow {\cal C}_0$. The only
possibility for the definition of $\hat f$ is obvious: any point
$P$ of $\hat X$ can be regarded as an IP and, thus, represented as
$I^-(c)$ for some future-directed chain $c$, so, $\hat f(P)$ is
defined as $I^-[f(c)]\in \hat Y$.
%Nevertheless, one
%has to check that $\hat f$ also preserves the chronology, and this
%depends on the definition of the chronology $\ll$ induced in the
%completions. As Harris chose $\ll = \ll^c$, his approach is
Unfortunately, $\hat f$ does not necessarily preserve the
chronology
(Fig. 6). %\ref{f6}).
So, one is forced to consider the subcategories ${\cal C}^{pd}$
and ${\cal C}^{pd}_0$ (of ${\cal C}$ and ${\cal C}_0$, resp.)
whose objects are past determined chronological sets. The
following categorically universal result can be then obtained
\cite{HaJMP}:
%By choosing $\ll = \ll^e$ this restriction
%seems not necessary.
%At any case, Harris' idea is clear:
\btheo
%Let ${\cal C}^*$ an {\em admissible} category  of chronological
%sets (either the  one ${\cal C}$ defined above or the  subcategory
%containing only past-determined chronological sets) with morphisms
%the future-continuous functions, and let ${\cal C}^*_0$ the
%subcategory of its complete objects.
In the subcategory of past-determined chronological sets ${\cal
C}^{pd}$: (i) each morphism $f: X\rightarrow Y$ extends naturally
and univocally to a unique morphism $\hat f: \hat X\rightarrow
\hat Y \in {\cal C}^{pd}_0$, (ii) future-completion $ \widehat \;
:{\cal C}^{pd}\rightarrow {\cal C}^{pd}_0$ is a functor, and (iii)
the standard future injection is a natural transformation between
functors. \etheo

The drawback of this universality is that it applies only in the
past-determined category. Nevertheless, Harris also defined a
natural  functor {\bf pd}: ${\cal C}\rightarrow {\cal C}^{pd}$ of
past-determination. A further study of this functor shows the
naturality of the GKP construction as a minimal  way of
``future-completing'' any past-regular, past-distinguishing
chronological set.

%That is, modulo the choice between $\ll^c, \ll^e$,  GKP
%construction is the minimal (universal, categorically unique) way
%of ``future-completing'' a past-regular, past-distinguishing
%chronological set.

%This yields a strong support to some parts of GKP construction:
%the definition of future preboundary points for any
%past-distinguishing spacetime  and the extended causal relations
%to $\hat M$, in the case of past-determined spacetimes\footnote{y
%dice que la eleccí\'on the $\ll^e$ es la buena, que' carajo.}.

\subsubsection{Chronological topology}

In order to define a topology in any $(X,\ll)$ associated to the
chronological relation, some first options must be disregarded:

\begin{remark}\label{reh}{\em
(a) As we have already commented in Subsection \ref{s3.1}, general
examples show that plain Alexandrov's topology yields
unsatisfactory properties in simple cases. A different type of
 example is cited by Harris \cite[Figure 2]{HaCQG} in the context
of future boundaries. This may be illuminating and will be
revisited below (Fig. 10). %\ref{f10}).
Consider  the chronological
space $X$ obtained by removing the negative spacelike semi-axis
$x\leq 0$ from $\LL^2$. Each removed point $(t=0,x)$ yields
naturally a boundary point $P_x (=I^-(0,x)) \in \hat
\partial X$. In the completion $\hat X$, the past of $(1,0)$ would
yield an Alexandrov open subset which contains $P_0$ but no other
$P_x$. So, $P_{1/n}\not\rightarrow P_0$. But, as we are looking
only at past sets, the convergence of the sequence seems
desirable.

(b) The GKP topology would not be appropriate, even though now the
problem of the identifications between future and past preboundary
points does not appear. Among the reasons discussed  previously,
the anomalous topology for $\hat\partial \LL^n$ (footnote
\ref{footh}) is especially unfavorable now, as $\LL^n$ is globally
hyperbolic and, thus, there are no identifications between future
and past boundary points. }\end{remark}
 \sm \sm
 Instead, Harris defined a {\em limit operator} $\hat L$: for any
sequence $\sigma =\{x_n\}_n \subset X$,

\be \label{limitop}
x\in \hat L(\sigma) \Leftrightarrow \left\{ \begin{array}{l} y\ll x \Rightarrow y\ll x_n %& \hbox{for large} \; n
\\
I^-(x) \varsubsetneq P (\in \hat M) \Rightarrow \exists z\in P: z
\not\ll x_n
\end{array}\right. \quad \hbox{for large }\; n$$
%(for large $n$)
or, equivalently, in Flores' reformulation:
$$x\in \hat L(\sigma) \Leftrightarrow \left\{ \begin{array}{l} I^-(x) \subset LI\{I^-(x_n)\}  \\
I^-(x) \, \hbox{is a maximal IP in } LS\{I^-(x_n)\}
\end{array}\right.\ee
where $LS, LI$ denote the lim-sup and lim-inf operators in set
theory, that is, for any sequence of subsets $A_n\subseteq X$,
$LS(\{A_n\})$ (resp. $LI(\{A_n\})$) contains any $x\in X$ which
belongs to infinitely many $A_n$ (resp. all $A_n$ for $n$
sufficiently large). This limit operator defines the closed
subsets for a topology, in fact:

\bdefi Let $X$ be a past-regular chronological set. The {\em
future-chronological topology} is the one such that: $C\subseteq
X$ is closed if and only if for any sequence $\sigma$ in $C$,
necessarily $\hat L(\sigma)\subseteq C$. \edefi

Notice: (i) as $X$ is chronologically separable, the so-obtained
topology is second countable (any second countable topology can be
characterized by such a limit operator), (ii) points are closed
because of past regularity, i.e., the future-chronological
topology is $T_1$.

The following properties suggest that  future-chronological
topology has been chosen properly\footnote{This topology depends
on $\ll$ and, so,  would depend  on the choice $\ll^c$ in
(\ref{eell}) for a completed spacetime. %(recall  Remark \ref{rh}).
Nevertheless, Harris emphasizes \cite{HaGeoD} that topological
properties will be independent of past determination --and, then,
on some subtleties of this choice.}:
\begin{enumerate}
\item For any future-directed chain $\sigma = \{x_n\}$: $x\in \hat
L(\sigma) \Leftrightarrow \{x_n\}\rightarrow x$.

\item The standard injection $X \hookrightarrow \hat X$ is a
homeomorphism into its image, and $X$ is dense in $\hat X$.

\item If $X$ is a strongly causal spacetime,  X has the manifold
topology and $\hat \partial X$ is closed in $\hat X$.

\item For $X=\LL^n$, the chr-boundary is the usual cone with its
natural topology.
\end{enumerate}
Additionally, the topology satisfies an interesting property of
quasi-compactness \cite[Sect. 5]{FH}: any sequence $\{x_n\}$ with
$LS(\{x_n\})\neq \emptyset$ admits a subsequence with a non-empty
limit.

Nevertheless, in spite of these nice properties, there is an
undesirable property. As we have already commented (Remark
\ref{reh}(a)), sets type $I^\pm(x)$ are not always open now. This
may be natural when only convergence of IP's  is being considered
(forgetting what happens for IF's). Nevertheless, this also
implies that, in the future chronology, points of the boundary are
not $T_2$ related to points of the spacetime (Fig. 10).
%\ref{f10}).

Recall also that there is no a general categorical result for the
future-chronological topology. The reason is the following.
Consider a future-continuous function $f:X\rightarrow Y$ (with $X,
Y$ past-determined) and its extension to the completions $\hat
f:\hat X\rightarrow \hat Y$. If $X, Y$ are endowed with  the
chronological topology, {\em the continuity of $f$ is not
sufficient to ensure the continuity of $\hat f$}. This happens so
easily (Fig. 7) %\ref{f7})
that seems to be an inherent obstruction to the categorical
approach: if the objects in the categories ${\cal C}, {\cal C}_0$
are endowed with the chr-topology, and the morphisms are assumed
to be also continuous, then future completion does not yield  a
functorial relation between ${\cal C}$ and $ {\cal C}_0$ (or
subcategories such as ${\cal C}^{pd}, {\cal C}^{pd}_0$). Moreover,
it is hard to think that some natural definition of the topology
of the c-boundary might satisfy such a universal property.

\begin{figure}[ht]\label{f7}
\begin{center}
\includegraphics[width=.8\textwidth]{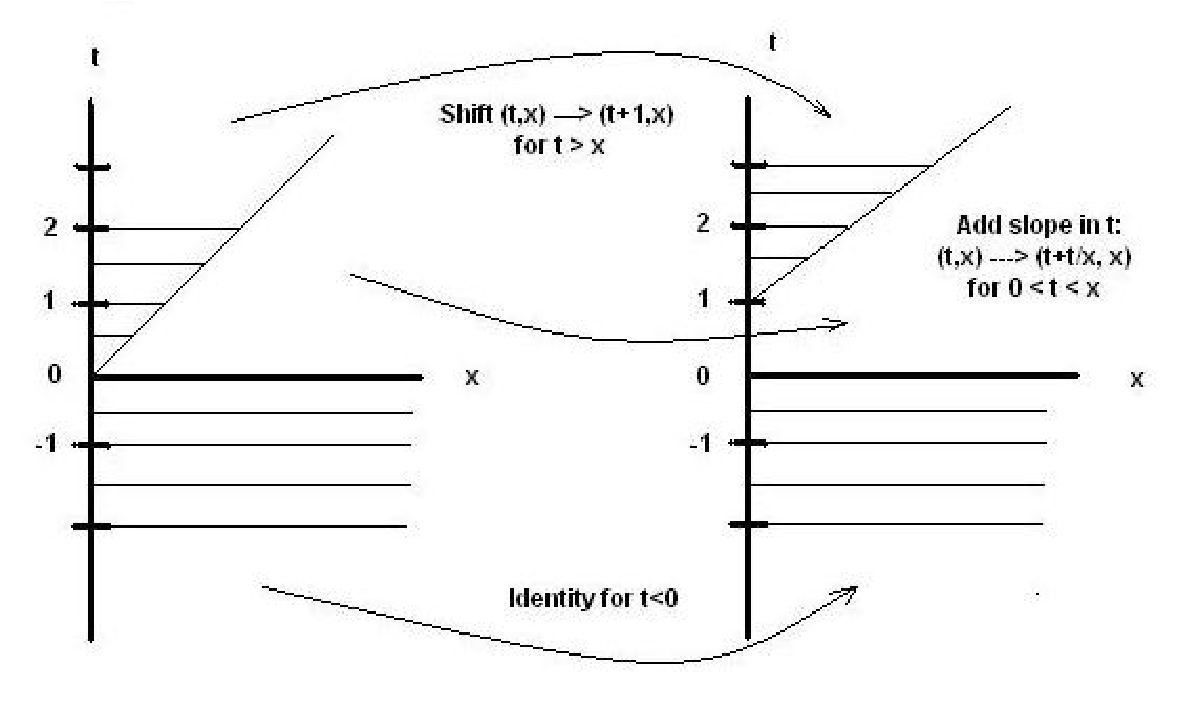}
\caption{Bijective chronological map of $\LL^1\times (0,\infty)
(\subset \LL^2)$ in itself. It is continuous with the chr-topology
(and future-continuous) but cannot be extended continuously to the
completions, as the origin ``expands'' to the vertical segment
$[0,1]\times \{0\}$ in the boundary, \cite[p. 569]{HaCQG}.}
%\medskip
%\parbox{5cm}{
%\footnotesize
%{\small Conformal embedding of Minkowski
% spacetime in the Einstein static universe.}

\end{center}
\end{figure}

At any case, Harris also proved that, in the subcategory of {\em
chronological sets with spacelike boundaries} (see Subsection
\ref{s5.5}), such a functorial relation is still possible, and the
universality of chr-boundary in the sense of categories, holds.
(Notice that in the case of spacelike boundaries, the problem of
identifications between future and past boundary points do not
appear.) Again, this yields a strong support for chr-topology.

 \subsection{Summary}

We retain the following elements of previous approaches: \ben
\item The GKP precompletion of $M^\sharp$ (spacetime + preboundary
points) is natural, but its chronology, topology and
identifications (either as in the GKP approach or in related
methods) become suspicious. \item
 For the future boundary $\hat
\partial M$, Harris universal properties for $(M,\ll )$ as a chronological set
holds. Even though there is some restriction (class of
past-determined chronological sets), perhaps some different
extension of the relation $\ll$ to the boundary might drop it. At
any case, the GKP construction of $\hat M$ finds a good support as
a point set and  a chronological set.
 \item
Harris limit operator $\hat L$ for the topology of $\hat M$
(avoiding GKP/Alexandrov topologies) becomes  natural as a
convergence of past sets. Nevertheless, as the future parts are
not taken into account, points of $\hat\partial M$ and $M$ may be
not $T_2$ related. Its universal properties are restricted to the
class of spacetimes with spacelike boundaries. However, this seems
inherent to the conformal character of the construction. Examples
suggest that this restriction would be unavoidable for the
universal properties of the topology, in any notion of the
c-boundary. \item At any case, it is necessary to relate $\hat
\partial M$ and $ \check \partial M$. Even though the procedure
was not clear, the Szabados  relation, which refined  Budic and
Sachs',
%$$
%P\sim_S F \Longleftrightarrow \left\{
%\begin{array}{l} F \quad
%\hbox{is included and  is maximal in} \quad \uparrow P \\
%P \quad \hbox{is included and is maximal in} \quad \downarrow F
%\end{array} \right. $$
became promising. \item Additionally, Harris \cite{HaNonlin}
introduced some tools to compute the boundary in standard static
spacetimes (refined later in \cite{FH, HaGeoD}) and quotient
spacetimes. \een And, of course, one also had a set of worrying
counterexamples, as those by Kuang and Liang!

\section{Intermission: holography and waves come into
play}\label{s2}

\subsection{AdS/CFT correspondence and boundaries}

%Now, let us recall very briefly some foundations on holography.
Very roughly, the idea of holography is that  physics in a region
is encoded by some fundamental degrees of freedom in the boundary
of this region --the ``hologram'' of the original region.  The
seminal idea appeared by studying  the entropy of black holes,
which is proportional to the area of the event horizon --a
surprising property, as entropy is an ``extensive'' magnitude and
one would expect proportionality to the volume. So, G. 't Hooft
\cite{tH} and L. Susskind \cite{Su} suggested that the nature of
quantum gravity might be holographic.
%Such a property
%extends to the information which can be extracted from the black
%hole and, if this were generalizable to any mass (not necessarily
%black holes) the volume will be at some extent an ``illusion''
%--its information would be codified as a hologram in a dimension
%less.
The most rigorous realization of the holographic principle is
Madacena's AdS/CFT correspondence. This is a conjectured
equivalence between: (i) string theory on a space (typically,
AdS$_5 \times \SSS^5$, or the product of anti de-Sitter by other
compact manifold), and (ii) a Quantum Field Theory (say, a
Conformal Field Theory) on the {\em conformal boundary} of this
space, which behaves as a hologram of inferior dimension.

We are interested in this holographic boundary. It is not
difficult to check that AdS$_5\times \SSS^5$ is conformally
equivalent to $(\R^6\setminus \{0\}) \times \LL^4$ and, thus, to
$\LL^{10}\backslash \LL^4$ (a timelike linear subspace $\LL^4$ is
removed from $\LL^{10}$). In fact, taking into account the
Poincar\'e representation of AdS$_5$, $g_P= (dy^2+g_{\LL^4})/y^2,
y >0$ the metric $g= g_P + g_{\SSS^5}$ in AdS$_5\times \SSS^5$ is
conformally  equivalent to:
$$
y^2 g  = dy^2+g_{\LL^4} + y^2 g_{\SSS^5} \equiv g_{\R^6}+g_{\LL^4}
\quad  \quad              \hbox{for } \; y>0.
$$
The conformal boundary obtained by means of the inclusion
(AdS$_5\times \SSS^5 \approx) \R^6\setminus \{0\} \times \LL^4
\hookrightarrow \LL^{10}$, restores the removed subspace $\LL^4$.
This is also the expected causal boundary and, so, it is not
relevant in this case which one of the two boundaries is chosen.

\subsection{Holography on plane waves}

Nevertheless, the situation will change when one consider
holography on plane waves. A priori, such a holography is
interesting because: (a) some plane waves provide exact
backgrounds for string theory (as all curvature invariants
vanish), and (b) every spacetime has a plane wave as a limit
(Penrose limit \cite{Pe}). But the string community was not truly
interested in this holography until  Berenstein, Maldacena and
Nastase detailed correspondence \cite{BMN} between 10d-string
theory on plane waves and 4d-super Yang-Mills theory. Moreover,
Blau, Figueroa-O'Farrill, Hull, Papadopoulos (BFHP) took a
lightlike geodesic in AdS$_5 \times \SSS^5$ which rotates on
$\SSS^5$, considered its Penrose limit and identified its dual in
field theory \cite{BFHP}.

Berenstein and Nastase \cite{BN} revised previous approach. They
found that, unexpectedly, the conformal boundary for BFHP wave is
1-dimensional. In their conformal embedding, the boundary was a
null line with a rather surprising role of ``future and past''
(Fig. 8). %\ref{f8}).
This additional reduction of the dimension
became relevant, as it opened new possibilities for the holography
(the holographic dual can be described by a quantum mechanical
system --a matrix model).

\begin{figure}[ht]\label{f8}
\includegraphics[width=.25\textwidth]{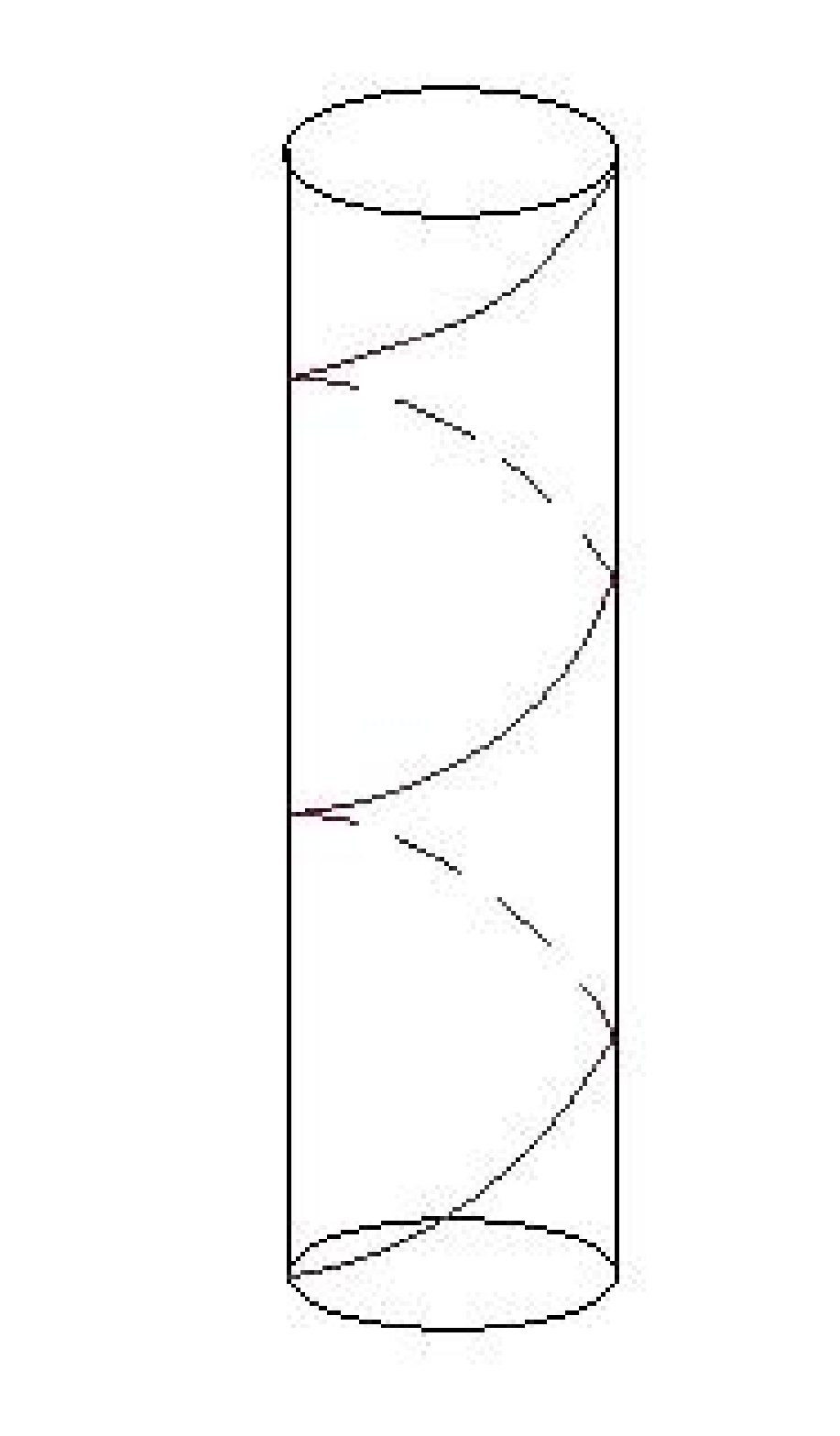}

%\medskip
%\parbox{5cm}{
%\footnotesize
%{\small Conformal embedding of Minkowski
% spacetime in the Einstein static universe.}
\caption{The depicted ``helix'' in Ein$_2$ represents the
1-dimensional boundary of BFHP wave. One can arrive at this
lightlike boundary by means of both, future and past directed
timelike curves.}
\end{figure}

\subsection{Is conformal boundary the right choice?}

Marolf and Ross \cite{MR1} realized that the usage of conformal
boundary was limited to very particular plane waves --say, the
conformally flat ones. As the causal boundary is conformally
invariant and applicable to many more waves, they proposed to
study such a boundary. Due to the results for the BFHP wave, the
problem of the identifications became essential and, so, they
introduced their own proposal, which will be discussed below.
Then, they proved that the 1-dimensional character of the
conformal boundary, is reproduced for the causal boundary not only
for  BFHP plane wave but also for other plane waves.

Even though they made a concrete choice of c-boundary, their
results held by assuming a minimal requirement, namely: for any
sequence $\{p_n\}_n\subset M$, if the PIP's $I^-(p_n)$ approach
(in any reasonable sense) some TIP $P$ and the PIF's $I^+(p_n)$
approach some TIF $F$ then $P,F$ represent the same point. This
puts forward the problem of the identifications.

As we will see next, Marolf and Ross study of c-boundary was
widely extended by Flores \cite{F}, and the systematic computation
of the c-boundaries of wave type spacetimes was carried out by
Flores and S\'anchez \cite{FS}. But, prior to this, the following
conclusion is clear:

\begin{quote}
In the AdS/CFT correspondence, the boundary of the spacetime was
used in an elementary way by means of conformal embeddings. But
when plane wave backgrounds come into play, this is no longer
satisfactory: the causal boundary must be used and the problem of
identifications becomes essential.
\end{quote}

\section{Reconstructing the c-boundary} \label{s3}

\subsection{Marolf-Ross seminal idea}

As pointed out in the Introduction, the key role of Marolf-Ross
 approach \cite{MR2} is to consider the c-boundary $\partial_{MR}
M$ as the set of all the pairs S-related (see (\ref{eSz}). That
is, for any strongly causal spacetime $M$ (or even any
distinguishing and regular chronological set) \be \label{mr}
\begin{array}{ll}
\partial_{MR}M = &\{(P,F): \,  P\; \hbox{is a TIP}, \;  F\; \hbox{is a TIF},
 \hbox{and} \; P\sim_S F\} \\ &
 \cup \{(P,\emptyset): \,  P\; \hbox{is a TIP},
 \hbox{and} \; P\not\sim_S F \; \hbox{for any TIF}\; F\}
 \\ & \cup \{(\emptyset, F): \,  F\; \hbox{is a TIF},
 \hbox{and} \; P\not\sim_S F \; \hbox{for any TIP}\; P
 \}
 \end{array}
\ee
 In particular, if $P\sim_S F \sim_S P'$ then $(P,F), (P',F)$ are
regarded as two different boundary points. Among the advantages of
this viewpoint, one is apparent: the chronological relation can be
extended in an obvious  way to the boundary (Fig. 9): %\ref{f9}):
\be \label{efp} (P,F) \overline{\ll} (P',F') \Leftrightarrow F\cap
P'\neq \emptyset , \quad \forall (P,F), (P',F')\in \partial_{MR}M
\ee

\begin{figure}[ht]\label{f9}
\includegraphics[width=.5\textwidth]{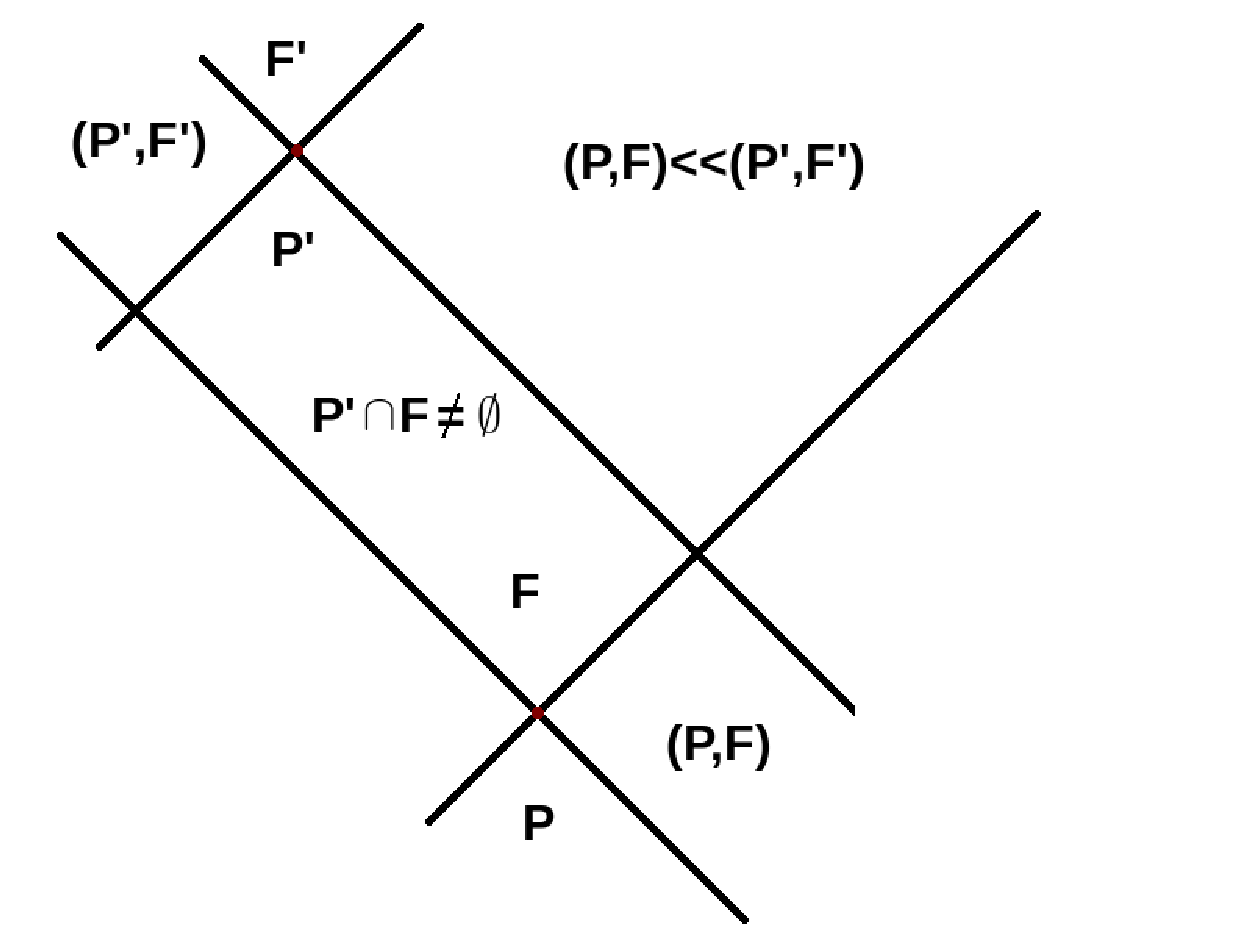}
\caption{When the points of the completion are regarded as pairs,
the chronological relation emerges naturally from (\ref{efp}).}
\end{figure}

This simple definition is satisfactory, as one can prove that
$\overline{\ll}$ satisfies: (i) is transitive, and (ii) does not
introduces additional (spurious) relations in $M$ --when applied
to PIP's and PIF's\footnote{Notice that such a $\overline{\ll}$
was essentially introduced by Szabados \cite{Sz1} but, under his
approach, conditions (i) and (ii) could not be fulfilled: the
problem was inherent in the fact that all the elements connected
by
$\sim_S$ represented the same boundary point, % i.e., $P\sim_S F
%\sim_S P' \Rightarrow P\equiv P' \equiv F$, yields this anomaly,
(see  \cite{MR2} or Example \ref{exfl} below).}.

\begin{remark}\label{rh} {\em The essential uniqueness of
the extended chronological relation for pairs contrasts with
previous approaches. For example, Harris' choice
%for the extended
%chronological relation
$\ll^c$ (see (\ref{eell})) became standard and natural, but we
have remarked the problems concerning past-determination. Taking
into account (\ref{efp}), a subtler extension $\ll^e$ would seem
also reasonable now:
%However, perhaps some alternatives might be of some interest, as
%the following:
%(1)  After the redefinition of the boundary (Section \ref{s3}),
%the following (a priori intrincated) completed chronological
%relation $\bar \ll^e$ becomes also natural:
$P\ll^e P'$ iff $P'$ is intersected by some IF, $F$, included and
maximal in $\uparrow P$, such that $P\sim_SF$. Apparently, this is
an extended chronological relation in $\hat M$ which does not
introduce spurious relations when $P, P'$ are TIPs.

Moreover, other possibilities of the chronological extension
(applicable to completions of chronological sets or  spacetimes)
might appear. For example, most of the striking differences
between the abstract chronological relations in $(X,\ll)$ and the
relation $\ll$ in a spacetime $M$, comes from the fact that the
latter is defined by means of paths. As any chronological space
can be also endowed with a  topology, one can generate a new
chronological relation which collects the idea of paths, either
directly (in spacetimes) or by mean of chains (in chronological
sets). So, start with some chronological set $(X,\ll)$
(eventually, take $X=\hat X, \ll = \ll^c$) and define: two points
$x, y \in X$ are {\em path-chronologically related}, $x \ll_p y$
if for every open covering $\{U_\alpha\}_\alpha$ of $X$ there
exists a chain $x=x_0 \ll x_1 \ll \dots \ll x_k=y$ such that each
consecutive pair $x_i, x_{i+1}$ is contained in some $U_\alpha$.
Any spacetime is {\em path-chronological} i.e. $\ll = \ll_p$; if
the extended relation $\ll^c$ in a topologized completion were
not, one could explore to change $\ll^c$ by $\ll^c_p$.

However such alternatives, which may be reasonable under previous
approaches, cannot be regarded as serious alternatives to
(\ref{efp}) for pairs.
%\footnote{De hecho, supongo que coincidirian
%con ella, cuando se formularan correctamente}

%Now, in a non-past determined chronological set, the local
%chronology $\ll_l^c$ associated to Harris $\ll^c$, may be also a
%natural chronology on the completion.}
}\end{remark}

So, Marolf-Ross' ideas introduce a new and exciting  viewpoint,
which was complemented by suggesting two possible topologies for
the boundary. But, in order to understand the exact role of MR
completion, let us consider first the general notion of completion
introduced by Flores \cite{F}.

\subsection{Flores' general completions: point set and chronological extension}

Let us retain the most general elements of Harris' setting,
extending when necessary past elements to future ones. That is,
consider, among other straightforward definitions, a {\em
chronological set} $(X,\ll)$, (future or past directed) {\em
chains} and (past and future) {\em limit operators} $\hat L,
\check L$ for any sequence $\sigma \subset X$ as in
(\ref{limitop}). For any future-directed chain $\eta=\{x_n\}_n$
the inclusions $I^-(x_n) \subset I^-(x_{n+1})$ and $I^+(x_n)
\supset I^+(x_{n+1})$ yield directly that the lim sup LS and lim
inf LI of each sequence $I^\pm (x_n)\}_n$ coincide and, then:
\bprop \label{ppio} Let $(X,\ll)$ be a (past and future) regular
chronological set. For any future (resp. past) chain $\eta$:
$$ \begin{array}{c}x \in \hat L(\eta)\Leftrightarrow I^-(x)=I^-(\eta) , \quad \quad
 x \in \check L(\eta) \Leftrightarrow I^+(x)\subset \uparrow \eta \;
 \hbox{{\rm and is maximal in }} \uparrow \eta
\\
 \hbox{{\rm (resp.}} \;\; x \in \check L(\eta)\Leftrightarrow
I^+(x)=I^+(\eta) , \quad \quad
 x \in \hat L(\eta) \Leftrightarrow I^-(x)\subset \downarrow \eta \;
 \hbox{{\rm and is maximal in }} \downarrow \eta  \hbox{{\rm ).}}
 \end{array}$$ \eprop

The natural notion of completion for $(X,\ll)$ must comprise that
any chain have an endpoint. So, we need two previous concepts: the
first one is the set-point space where this completion makes
sense,  the second is the notion of {\em endpoint of a chain}.

Let $X_p, X_f$ be, resp., the sets of all the   past sets
($P\subset X$, $P = I^-(P)$) and future sets ($F\subset X$, $F =
I^+(F)$) --not necessarily IPs or IFs. Assuming that the
chronological set is {\em weakly distinguishing} (i.e., either
future or past distinguishing) the original chronological set
injects naturally,
$$
\mathbf{i}: X\rightarrow X_p\times X_f , \quad \quad x
\rightarrowtail (I^-(x), I^+(x)).
$$
So, as a first conclusion,  any completion $\bar X$ would satisfy
$\mathbf{i}(X) \subset \bar X \subset X_p\times X_f$ as a point
set.

The notion of endpoint is subtler. Recall that we have not defined
a topology yet, even though the operator $\hat L$ does define a
natural topology in $\hat X$. In principle, an endpoint is not the
same thing that a {\em limit}, but Proposition \ref{ppio} is a
clear guide. The essence of an endpoint $(P,F)$ for, say, a
future-directed chain $\eta$ deals with future convergence, and
naturally one must impose $P=I^-(\eta)$. Nevertheless, some
control is necessary for the $F$ component, and this is carried
out by means of the operator $\check L$. As this operator yields
IFs, let us introduce first a {\em decomposition operator}, dec.
By using Zorn's lemma, it is easy to show that any $P\in X_p$ can
be written as an union $P=\cup_\alpha P_\alpha$, where the set
$\{P_\alpha: \alpha\in I\}$ contains all the IP's included in $P$
which are {\em maximal} under the relation of inclusion. Consider
the following operator, which applies on past sets and, dually, on
future sets: \be \label{dec} \hbox{dec}(P) = \{P_\alpha: P_\alpha
\, \hbox{is a maximal IP included in} \, P , \forall \alpha\in I
\}, \quad \quad \hbox{for any} \, P\in X_p .\ee Now, we have the
elements for the definitions.

\bdefi \label{d5.2}Let $X$ be a weakly distinguishing
chronological set $X$.

(a) A pair $(P,F)\in X_p\times X_f$ is an {\em endpoint} of a
future (resp. past) chain $\eta \subset X$ if
$$
P=I^-(\eta ), \; \hbox{dec} (F) \subset \check L (\eta) \quad
\quad (\hbox{resp.} \; \; F=I^+(\eta ), \; \hbox{dec} (P) \subset
\hat L (\eta )).
$$

(b) $X$ is {\em complete} if any chain $\eta$ in $X$ has some
endpoint
in $X$. \ %The required definition is then:

(c) A set $\bar X, \mathbf{i}(X)\subset \bar X \subset X_p\times
X_f$ is a {\em completion} of $X$ if any chain $\eta$ in $X$ has
some endpoint in $\bar X$. In this case, the {\em extended
chronological relation} $\overline{\ll}$ on $\bar X$ is defined as
in (\ref{efp}). \edefi

The boundary of the completion is then $\partial X:= \bar
X\backslash \mathbf{i}(X)$. Some natural properties show the
consistency of the definitions (see \cite{F} for detailed proofs):
(i) $(\bar X, \overline{\ll})$ becomes a weakly distinguishing
chronological set, (ii) the injection $\mathbf{i}$ is a
chronological map between $(X,\ll)$ and $(\bar X,
\overline{\ll})$, (iii) the image $\mathbf{i}(X)$ is
chronologically dense in $\bar X$, (iv) no spurious chronological
relations are introduced by $\overline{\ll}$ in $\mathbf{i}(X)$,
and (v) any completion $(\bar X, \overline{\ll})$ of $X$ is a
complete chronological set. We emphasize the natural role of
$\overline{\ll}$ for these properties. For example, by using a
different extended chronology, Harris had constructed a spacetime,
which suggested a generic counterexample to property
(v) --see footnote \ref{foot}.% \cite[Appendix]{HaJMP}.

\subsection{Minimal completions and Marolf-Ross one}

Previous notion of completion is very general and includes not
only Marolf-Ross one but, for example, GKP precompletion (looking
each IP $P$ as the pair $(P,\emptyset)$ and each IF $F$ as
$(\emptyset, F)$). In order to have more efficient completions, we
must restrict to appropriate {\em minimal} ones.

Consider the set ${\cal C}_X$ of all the completions of a weakly
distinguishing $X$. Now, delete from ${\cal C}_X$ those
completions which are still a completion when some of the points
of the boundary is removed. The resulting  subset ${\cal C}^*_X$,
is not empty (for example, GKP precompletion remains there), and a
partial ordering $\leq$ is defined as follows. Let $\bar X^I, \bar
X^J \in {\cal C}^*_X$, we put $$\bar X^I \leq \bar X^J$$ if for
each $(P_i,F_i)\in \partial X^I$ there exists some $S_i \subset
\partial X^J$ such that the following properties are fulfilled:
\begin{itemize}
\item[(a)] The set of all the $S_i$'s is a partition of $\partial
X^J$.

\item[(b)] If some chain $\eta$ in $X$ has an endpoint in $S_i$
then $(P_i,F_i)$ is also an endpoint of $\eta$.

\item[(c)] If some $S_i$ contains only a point
%\footnote{De esta
%propiedad y de la anterior se deduciria que los dos puntos son
%endpoints de las mismas cadenas, sin necesidad de imponerlo ahi',
%no?}
$S_i=\{(P,F)\}$ then dec$(P_i)\subset $ dec$(P)$ and
dec$(F_i)\subset $ dec$(F)$.\eit Recall that (a) implies the
existence of an onto map $\Pi: \partial  X^J\rightarrow
\partial  X^I$. By (b) if all the boundary points in $S_i$ were
 replaced by $\Pi(S_i)=(P_i, F_i)$, then one would still obtain  a
completion (however, recall that $S_i$ has at most two points).
Finally, (c) ensures that, if $(P,F)$ is replaced by $(P_i,F_i)$
then the latter is simpler than the former (and the chains will
have still endpoints). The existence of minimal elements for the
partial order $\leq$ is guaranteed by using Zorn's lemma.

Such  minimal completions are also called {\em chronological
completions} in  \cite{F}. They fulfill very satisfactory
properties
%in the case that $X$ is distinguishing and regular or,
for  strongly causal spacetimes (see \cite[Th. 7.4]{F}):

\btheo \label{tlamadre} Let $M$ be a strongly causal spacetime. A
completion $\bar M \subset M_p \times M_f$ is minimal if and only
if its boundary $\partial M$ satisfies the following properties:

(i) Every TIP and TIF in $M$ is the component of some pair in
$\partial M$.

(ii) If $ (P,F) \in \partial M$ and $P\neq \emptyset \neq F$, then
$P$ is a TIP, $F$ is a TIF, and $P\sim_S F$.

(iii) If $ (P,F) \in \partial M$ and $F=\emptyset$ (resp.
$P=\emptyset$) then $P$ (resp. $F$) is a (non-empty) terminal
indecomposable set and is not S-related to any other set.

(iv) If $(P,F_1), (P,F_2) \in \partial M$ and $F_1\neq F_2$
(resp., $(P_1,F), (P_2,F) \in \partial M$ and $P_1\neq P_2$) then
$F_i$ (resp. $P_i$), $i=1,2$, does not appear in another pair of
$\partial M$.

Moreover, Marolf-Ross completion (\ref{mr}) is the maximum
completion (the biggest as a subset of $M_p\times M_f$) which
satisfies properties (i), (ii) and (iii).

\etheo

Notice that property (i) is necessary to make $\bar M$ a
completion, and (iv) to make it minimal. The appearance of
properties (ii) and (iii) become specially relevant. First,
%\begin{quote}
 {\em Szabados relation appears naturally, it is not
imposed a priori as in previous approaches.}
%\end{quote}
And, second, as a straightforward consequence of Marolf-Ross'
definitions: \bcoro Let $M$ be a strongly causal spacetime, $\bar
M$ a minimal completion and $\bar M_{MR}$ Marolf-Ross one. Then
$\bar M \subset \bar M_{MR}$, and the inclusion is a chronological
map\footnote{In fact, it is also future (and past) continuous
(Defn. \ref{dmorf}), and a chronological isomorphism onto its
image.}. \ecoro

\brema\label{rmr} {\em (a) The appearance of S-relation $\sim_S$
for all the points in the completion (Theorem \ref{tlamadre}), is
a more general fact for (say, regular, distinguishing)
chronological sets. Nevertheless, there is an important reason for
the restriction to the class of strongly causal spacetimes.
According to \cite[Proposition 5.1]{Sz1}: {\em (i) A spacetime $M$
is strongly causal iff $I^+(p)\sim_S I^-(p)$ for all $p\in M$, and
(ii) if $M$ is strongly causal then $I^-(p) \sim_S F$ iff
$F=I^+(p)$} (and analogously for $I^+(p)$). That is, when $M$ is
not strongly causal, then $(I^-(p),I^+(p))$ is not always a pair
taken into account in Marolf-Ross completion (such pairs should be
added directly in order to get a completion), and pairs type
$(P,I^+(p))$, $(I^-(p), F)$ with $p\in M$ and $P, F$  terminal
sets, may appear.

(b) For a strongly causal spacetime $M$, one would be tempted to
consider $\bar M_{MR}$ as the smallest completion which includes
all the minimal ones. Unfortunately, this property do not hold:
even when there exist only one minimal completion,  $\bar
M_{min}$, Marolf-Ross one may be strictly bigger than $\bar
M_{min}$ (see Example \ref{exfl}(2) below).

(c) However, it is clear that properties (i), (ii) (iii), which
appeared axiomatically in previous approaches, now receive a
strong support. So,  call  any completion which fulfills the three
properties an {\em  admissible completion}\footnote{Such a
completion, as well as some results below, will be developed in
\cite{FHS}.}. Then, Marolf-Ross completion is singled out as the
maximum admissible completion of $M$.}\erema

The following examples, taken essentially from \cite[Example
9]{F}, \cite[Appendix]{MR2}, stress the role of Marolf-Ross
completion as a non-minimal one, as well as the importance of
pairings for the consistency of $\overline{\ll}$.

\bexam \label{exfl}{\em (1) Delete from $\LL^3$ the   coordinate
semi-planes: $XT^+=\{(t,x,y=0): 0\leq t \}$ and $YT^-=\{(t,x=0,y):
 t \leq 0\}$. Clearly,  there are four terminal sets associated
with the removed origin: two TIPs $P_1, P_2$ (corresponding to the
$x>0$ and $x<0$ directions, resp.)  and two TIFs $F_1, F_2$
(corresponding to  $y>0$ and $y<0$, resp.); moreover $P_i\sim_S
F_j$ for $i,j=1,2$. Marolf-Ross completion includes then the four
pairs $(P_i,F_j)$ as boundary points. Nevertheless, minimal
completions are obtained by attaching only $(P_1,F_1), (P_2,F_2)$
or only $(P_1,F_2), (P_2,F_1)$. Both are naturally included in
Marolf-Ross. In this concrete example, the two minimal completions
are isomorphic in a natural sense, but this is not expected in
general.
%\footnote{Al final se puede poner un ejemplo  \cite{FHSst}}.

(2) Remove in previous example $X^+Y^+=\{(t=0,x,y): 0\leq x, 0\leq
y\}$. Now, $P_1\not\sim F_1$ and Marolf-Ross retains the other
three pairs. Nevertheless,  $(P_1,F_2), (P_2,F_1)$ yield a minimal
completion, which is {\em unique and strictly smaller} than
Marolf-Ross one.

(3) Modify slightly example (1) by removing also $Y^+T^+
=\{(t,x=0,y): 0\leq y, 0\leq t\}$. Now, $F_1$ splits naturally in
two: $F_1'$ which contains $p=(2,1,1,)$ and $F_1''$ (thus
$(2,-1,1)\in F_2''$). The past set $P_2$, which contains $-p$, is
not directly S-related with $F_1'$ according to (\ref{eSz}). In
fact, no point with $t,x,y <0$ is chronologically related to any
with $t,x,y>0$. None of the pairs would be chronologically related
according to $\overline{\ll}$. Nevertheless, according to Szabados
identifications, all the terminal sets $P_1, P_2$, $F_2$, $F_1'$,
$F_1''$ would collapse to a single boundary point $Q$. Therefore,
according to Szabados, one would have $-p\overline{\ll} Q
\overline{\ll} p$. So, either a spurious relation would be
introduced in the spacetime, or the transitivity of
$\overline{\ll}$ would be violated. No completion in the sense of
Defn. \ref{d5.2}  have such an important drawback.
%(By the way, Marolf-Ross completion is not minimal here.)
}\eexam

\subsection{The topology}

Up to now, the development of the c-boundary has depended on three
basic ideas: (a) completions as subsets of $X_p\times X_f$, (b)
the extended chronological relation (\ref{efp}), and (c) Defn.
\ref{d5.2} of endpoint. The first one was just a general setting
for completions, which includes any previous one. The second one,
even though non-trivial, is the natural and apparently unavoidable
choice (transitive and free of spurious relations) inside this
setting. Only the third one could be thought as a ``choice'' among
some conceivable alternatives. Nevertheless, Defn. \ref{d5.2}
allows one to have endpoints in such a  general way, that any
alternative definition would add additional requirements  in order
to have an endpoint of a chain --that is, the corresponding
completion would be a particular case of those already defined.
The other two relevant definitions --minimal and Marolf-Ross
completions-- are selected just by its desirable properties and
uniqueness.

For the definition of the topology there exists, in principle, a
bigger arbitrariness of the choice. However, we will see that
Flores' choice has no reasonable alternative from three
viewpoints: (A) the a priori natural choice in the setting, (B)
the good mathematical properties obtained a posteriori, and (C)
the uniqueness properties. But, prior going further, let us start
with an example of topological choice.

This example was already introduced in Remark \ref{reh}(a) (see
Fig. 10). %\ref{f10}).
Here, one has, for instance, $P\in \hat L(\{P_{-1/n}\})$, as we
focus only on the future boundary $\hat\partial M$. But this also
suggests that $Q=(P,F)\in \partial M$ must lie in the limit of the
 sequences $\{x_n\}, \{Q_n\}$. Marolf and Ross \cite{MR2} gave
two topologies. For the preferred one (to be described below),
$\{x_n\}, \{Q_n\} \not\rightarrow Q$, but this was regarded as a
non-desirable property. So, they introduced an alternative coarser
topology and, then, both sequences converge\footnote{Really,
Marolf and Ross \cite[Figure 5]{MR2} considered a more
sophisticated example which comes from \cite{GKP}. Nevertheless,
the intuition for the convergence of the sequences is similar to
the simplified example above -and $I^\pm(z)$ is not open for some
points of their alternative topology.}. Flores \cite[Example 10.4
and p. 631]{F} claims that such convergences are not natural, as
the intuition of convergence does not take into account what is
happening with the future parts of the pairs. Now, recall that
there exists a powerful reason to support this last opinion: in
the completion $I^-(z) (\equiv I^-(z,\bar M))$ contains $Q$ and
none of the points in the sequences. So, if $I^\pm(z)$ must be
open (and this is truly a natural requirement for the topology!)
then the sequences cannot converge.

\begin{figure}[ht]\label{f10}
\begin{center}
\includegraphics[width=.8\textwidth]{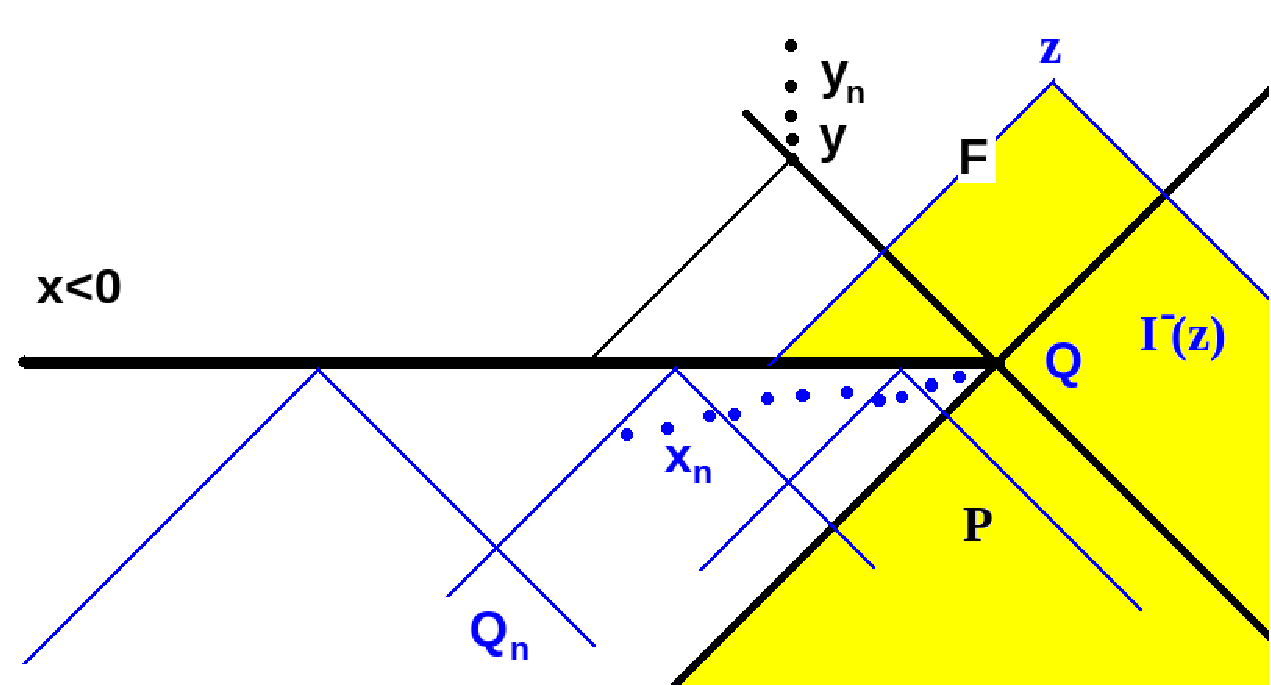}
\caption{In $\LL^2 \backslash\{x\leq 0\}$ each point of the
removed semi-axis  yields two boundary points $(P_x,\emptyset),
(\emptyset,F_x)$, when $x<0$, and a unique point $Q=(P,F)$ for
$x=0$. This spacetime $M$ exemplifies remarkable properties: (a)
$I^-(y)$ and $P$ lie in $\hat L(\{y_n\}_n)$, that is, $y\in M$ and
$P\in \hat\partial M$  are not $T_2$ related for the future
chronological topology, (b) if $I^-(z)$ is open in $\bar M$,
neither the sequence $\{Q_n=(P_{-1/n},\emptyset )\subset
\partial M\}_n$ nor the depicted one $\{x_n\}_n \subset M$  can converge to
$Q$, (c) the lack of global hyperbolicity of $M$ appears just
because the boundary is timelike at $Q$.}
%\medskip
%\parbox{5cm}{
%\footnotesize
%{\small Conformal embedding of Minkowski
% spacetime in the Einstein static universe.}

\end{center}
\end{figure}

Now, we describe Flores' choice. Along the discussion we will
focus on $(X,\ll)$ which are both {\em regular and
distinguishing}, even though some items do not require these
(mild) conditions.

\sm

\noindent {\em (A) The natural notion of convergence and
topology.} Recalling Defn. \ref{d5.2} and Prop. \ref{ppio}, $x$ is
an endpoint of a  (past or future-directed) chain $\sigma$ iff
\be\label{etop} \hbox{dec}(I^-(x))\subset \hat L(\sigma), \quad
\quad \hbox{dec}(I^+(x))\subset \check L(\sigma).\ee A minimum
requirement for the topology is that the endpoint $x$ will be also
the limit of the chain $\sigma$. So, the simplest natural choice
which fulfills this requirement is, obviously:

\bdefi \label{dtop} Let $(X,\ll)$ be regular and distinguishing.

(1) Given a sequence  $\sigma$ in X, we say  that {\em $x$ is a
limit of $\sigma$, or belongs to the limit operator $L$ on
$\sigma$} (denoted $x\in L(\sigma)$), if both inclusions
(\ref{etop}) hold.

(2) A subset $C\subset X$ is {\em closed} if it contains all the
limits of all the sequences contained in $C$.   The {\em
chronological topology} on $(X,\ll)$ is the one generated by the
closed sets. \edefi

Of course, this topology can be defined in any chronological set,
%(not necessarily regular or distinguishing),
as in \cite{F}. Perhaps, some alternative option would be
conceivable if $X$ were not regular, as in this case endpoint and
limit do not coincide for chains. But in the regular case, Defn.
\ref{dtop} implies directly that both concepts coincide, and any
alternative would be more complicated. \sm

\noindent {\em (B) Good properties of the topology.} The
properties fulfilled by the chronological topology were studied
systematically in \cite{F}. Some of them are valid for any
completion, but, as one could expect, the full desirable
properties hold only for the minimal completions of strongly
causal spacetimes  or, with some more generality, for their
admissible completions (see Remark \ref{rmr}(c)). Summing up:

1.- For any completion $\bar X$: (i) $\mathbf{i}: X \rightarrow
\bar X$ is a topological imbedding, and $\mathbf{i}(X)$ is dense
in $\bar X$, (ii) any chain in a complete chronological space (in
particular,  in any completion) has a limit, and (iii) for a
strongly causal spacetime $M$, the topology of the spacetime
coincides with the chronological topology, and any inextensible
timelike curve in $M$ has a limit in $\partial M$.

2.- For any admissible completion $\bar M$ of a strongly causal
spacetime $M$, including minimal and Marolf-Ross ones: (i) the
boundary $\partial M$ is a closed subset of $\bar M$, (ii) $\bar
M$ is $T_1$, (iii) if two points $Q, R\in \bar M$ are
non-Hausdorff related then both lie in $\partial M$, and (iv)
$I^\pm(Q)$ is open for any $Q\in \bar M$ (for this last property,
see also the forthcoming study \cite{FHS})\footnote{Notice that
property (ii) is not fulfilled by Marolf-Ross', (iii) is not
fulfilled by Harris', and (iv) is not fulfilled by Harris' and the
alternative Marolf-Ross topology.}.

\sm

\noindent {\em (C) Uniqueness properties.} The admissible
alternatives to the topological chronology for $(X,\ll)$ will be
studied in detail in \cite{FHS}. Here, we announce only that any
admissible topology must satisfy both, compatibility with $\ll$
(concretely, $I^\pm(x)$ must be open and the endpoints of chains
must be also limits of the chains) and compatibility with the
set-point limit operators (consider a converging sequence
$\{Q_n\}_n\rightarrow Q$,  looked in $X_p\times X_f$ i.e. $\{(P_n,
F_n)\}_n \rightarrow (P,F)$; if it happened $P\subset P'\subset
LI(P_n), F\subset F'\subset LI(F_n)$, then any neighborhood of
$(P,F)$ should contain $(P',F')$). Among such topologies, the
chronological one is selected as the one with best properties
%\footnote{En fin, ser las ``coarsest'' podria interpretarse
%tambien como una ``buena property'' supuesto que no se pueda
%deducir de la convergencia de las cadenas a sus puntos finales (o
%incluso de la equivalencia de estos dos conceptos en
%espaciotiempos fuert. causales - o e.chron. regulares y
%distinguidores)}
(as separation $T_1$).

In particular, it is also interesting to consider the role of the
main Marolf-Ross topology  in this setting. So, let $M$ be now a
strongly causal spacetime and $\bar M$ its Marolf-Ross completion
(or any admissible one). For any subset $S\subset \bar M$ consider
the following subsets of $\bar M$:

$$\begin{array}{l}
L^+_{IF}(S)=\{(P,F)\in \bar V: F\neq \emptyset, F\subset
\cup_{(P',F')\in S}F'\}, \\ L^-_{IP}(S)=\{(P,F)\in \bar V: P\neq
\emptyset, P\subset \cup_{(P',F')\in S}P'\}\end{array}
$$
If $L^+_{IF}(S), L^-_{IP}(S)$ are computed in $M$ they yield the
closures of $I^\pm(S)$. Nevertheless, the restriction $F\neq
\emptyset$ in $L^+_{IF}(S)$ (which is clearly necessary, as the
empty set is included in the subsequent union of subsets) makes it
so that no pair $(P,\emptyset)$ belongs to $L^+_{IF}(S)$. This
situation is rectified by introducing two operators $Cl_{FB},
Cl_{PB}$, the closures in the  future and past boundaries:
$$\begin{array}{l}
Cl_{FB}(S)=S\cup \{(P,\emptyset)\in \bar M: P\in \hat L(P_n) \;
\hbox{for some sequence} \; \{(P_n,F_n)\}_n\subset S\},
\\Cl_{PB}(S)=S\cup \{(\emptyset,F)\in \bar M: F\in \check L(F_n) \;
\hbox{for some sequence} \; \{(P_n,F_n)\}_n\subset S\}.
\end{array}
$$
Then, the Marolf-Ross topology is generated by using as closed
subsets: %\footnote{Realmente bastaria con poner
%$Cl_{FB}(L^+_{IF}(S))$ porque $S$ está incluido o en $L^+_{IF}(S)$
%o en $Cl_{FB}(L^+_{IF}(S))$ y, si uno toma
%$Cl_{FB}(Cl_{FB}(L^+_{IF}(S)))$ no obtiene nada nuevo, no? Tal vez
%conviniera poner esta definicion mas simplificada en \cite{FHS} y
%justificarla, para dejar claro que nuestro estudio va en serio}:
$$
L^+(S)=Cl_{FB}(S\cup L^+_{IF}(S)), \quad L^-(S)=Cl_{PB}(S\cup
L^-_{IP}(S)).
$$
Even though this choice of topology is natural, the possibility of
different alternatives, especially in the definition of the
operators $Cl_{FB}, Cl_{PB}$, was already pointed out by Marolf
and Ross. It is also worth pointing out that the  definition of
these operators in \cite{MR2} does not use our limit operators
$\hat L, \check L$, even though it becomes equivalent for pairs
$(P,\emptyset), (\emptyset, F)$ (see \cite{FHS}). Very remarkably,
the unique differences of limits of sequences between Marolf-Ross
and Flores topologies, may occur only when one of the components
$(P,F)$ of the candidate to limit is empty;  in this case,
Marolf-Ross topology may be non-$T_1$ (a typical example would by
Fig. 5); %\ref{f5});
this will be proven in a further study
\cite{FHS}. So, as commented above,  Flores' choice fix the best
behaved topology, among the admissible ones in the setting of
completions.

\subsection{Note: timelike boundary and extended causal relation}\label{s5.5}

 The extension $\bar{\ll}$ of $\ll$ to the completion $\bar M$
 is  not only natural for the construction of the c-boundary, but
 also a source of information on the original chronological set.
For example, a strongly causal spacetime $M$, or even just a
causal one, is {\em not} globally hyperbolic if and only if there
exists a boundary point $Q=(P,F)\in \partial M$ which satisfies
$x\overline{\ll} Q \overline{\ll} y$ for some $x,y\in M$ (see
\cite[Th. 9.1]{F}, and recall \cite{BeSa2}, \cite{SaERE},
\cite[Th. 5.11]{FH}). This condition is equivalent to $F\neq
\emptyset \neq P$, and  such pairs constitute the {\em timelike
boundary} $\partial_0M$. As commented in the Introduction,
$\partial M$ splits then in three subsets, $\partial_0M$,
$\partial_+M$ (composed by pairs $(P,\emptyset)$) and
$\partial_-M$ (pairs $(\emptyset, F)$). In the globally hyperbolic
case, the spacetime splits smoothly as an orthogonal product $\R
\times S$ (see \cite{BeSa}), and the boundary is expected to
reflect asymptotic directions, event horizons, and other causal
elements. Notice that even in this case ($\partial_0M=\emptyset$)
the parts $\partial_\pm M$ may look like very different; for
example, in a standard half-cylinder $(\R^-\times
\SSS^1,´-dt^2+d\theta^2)$, $\partial_-M$ is a point and
$\partial_+M$ a circle $\SSS^1$.

Up to now, we have considered only the chronological relation
$\ll$ and its extension to the boundary $\overline{\ll}$. However,
it is also natural to wonder for an extended causal relation and,
then for the ``lightlike'' parts of $\partial_\pm M$. Following
\cite{GKP}, there is a way to construct a causal relation
$\leq^{\ll}$ from any chronological relation $\ll$, namely:
$$ x\leq^{\ll} y \Leftrightarrow I^-(x)\subset I^-(y) \;
\hbox{and} \; I^+(x)\supset I^+(y) .$$ For spacetimes,
$\leq^{\ll}$ agrees $\leq$ in causally simple spacetimes (but not
in $\LL^2\backslash \{0\}$, consider $(-1,-1)$ and $(1,1)$). So,
if we consider the completion $\bar M$ of a spacetime $M$ and take
the associated causal relation $\leq^{\overline{\ll}}$, this will
introduce spurious causal relations in $M$. Harris \cite{HaCQG}
also introduced a notion of {\em chronological set $X$ with only
spacelike boundaries}. In the regular case, this comprises two
conditions: (a) if $P\in \partial\hat X$ then $I^-(P)$ (computed
with $\ll^c$ in (\ref{eell})) is not included in $I^-(Q)$ for any
$Q\in \hat X$, and (b) $\partial X$ is closed (for Harris
topology) in $\hat X$. But notice that condition (b) is technical,
and condition (a) (which surely would mean $P \not\leq^c Q$ for
any sensible definition of a extended causal relation $\leq^c$) is
very restrictive. So, it does not suggest how to introduce a
causal relation in $\bar M$. Finally, Marolf and Ross studied some
possible alternatives, and the associated problems regarding
transitivity and reflexivity \cite[Sect. 3.2]{MR2}; one could add
even some more alternatives
(recall Remark \ref{rh}). %\footnote{realmente, no he pensado  hasta
% que' punto las nociones de  {\em local chronology} en Remark
%\ref{rh} o su extension natural  {\em local causality} podrian
%valer para algo}.
However, as Marolf-Ross pointed out, this problem is not essential
at this level. In the case of waves with a $1$-dimensional
boundary, some properties suggest that this boundary might be
regarded as lightlike. But one can postpone a proper definition
until new issues make it necessary.

\section{Computation of the boundary of the waves}\label{s4}

Next, we will sketch how to compute the c-boundary in the case of
wave-type spacetimes. The computation of the boundary for the
simplest case, i.e., product spacetimes $\LL^1\times S$
is not trivial \cite{AF}. %In this case, the computation of the boundary is
%simplified, as each IP is  S-related with at most an IF --there is
%a unique minimal completion, equal to Marolf-Ross one.
 Busemann functions are required \cite{HaNonlin}, and subtleties at the topological level appear \cite{FH} (see also
 \cite{FHSst}). However, the
result is very intuitive and, at least at the level of the
boundary as a point set, more or less expected. Nevertheless, as
Marolf and Ross pointed out, plane waves are a physical example
were the necessity of a rigorous definition of the c-boundary is
stressed --in particular, the unique guidance to make pairs
$(P,F)$ is the abstract Szabados relation. However, in the case of
reasonably physical wave-type spacetimes, some of the subtleties
of this boundary do not appear. In fact, Marolf-Ross' is also a
minimal completion and, thus, the unique admissible one.

\subsection{The class of wave-type spacetimes}

Next, consider the class of spacetimes, namely, Mp-waves:

$$\m = \mo \times \R^2 , \; \;
\langle\cdot,\cdot\rangle_L = \langle\cdot,\cdot\rangle_0 -F(x,u)\
du^2 - 2\ du\ dv
$$
where $ (\mo, \langle\cdot,\cdot\rangle_0)$ is a $n_0$-dimensional
Riemannian manifold, $ (v,u)$ are the natural coordinates of
$\R^2$ and $F : \mo \times \R \to \R $ is a smooth function.
Particular cases are the {\em pp-waves} (plane fronted waves with
parallel rays), where $(\mo,\langle\cdot,\cdot\rangle_0)$ is just
$\R^2$ (or, in general, we also consider $\R^{n_0-2}$); the vacuum
or {\em gravitational} pp-waves correspond to $\Delta_xF(x,u)
\equiv 0$. Especially, plane waves are characterized as those with
$F(x,u)$ a quadratic form in $x$ for each $u$. For more background
on these spacetimes, see \cite{CFS, FScqg, FShonnef}.
%(x_1,x_2)
%\left(
%\begin{array}{rr}
%f_1(u) & g(u) \\
%g(u) & -f_2(u)
%\end{array}
%\right) \left(
%\begin{array}{l}
%x_1 \\
%x_2
%\end{array}\right)
%$ ($f_1\equiv f_2$: gravitational case).
%; $f_1, f_2, g$ compact
%support: sandwich wave)

The progress in the computation of the c-boundary  is contained
in: \bit \item  Berenstein and Nastase's seminal work \cite{BN} on
the conformal boundary for plane waves, under the assumptions: (i)
locally symmetric ($F(x,u) \equiv F(x)$) and (ii) conformally flat
(the eigenvalues of $F$ are positive and equal).

\item Marolf and Ross article \cite{MR1}, under their causal
boundary approach (see also \cite{MR3}). Again, only plane waves
are considered, but the assumption on conformal flatness is
removed.

\item Works by Hubeny, Rangamani and Ross \cite{HR,  HRRa, HRR}:
TIP's and TIF's are computed for some specific pp-wave backgrounds
for strings.

\item Flores and the author's   general framework for the
c-boundary in any wave-type spacetime \cite{FS}. This includes
relevant cases of pp-waves, and will be described below. \eit

\subsection{Previous remarks}

For the general approach in \cite{FS}, notice as previous
considerations:

(i) The  difficult part will be to compute TIPs, TIFs and their
common future or pasts, $P,F,\uparrow P, \downarrow F$. In fact,
 each $P$ will be $S$-related with at most one $F$ in the interesting
 cases.

(ii) In principle, one must compute $I^\pm(\gamma)$, $\uparrow
I^-(\gamma)$ $\downarrow I^+(\gamma)$ for any timelike  $\gamma$.
Nevertheless,  it will be enough to consider lightlike
curves\footnote{However, a subtlety appears: in general, the sets
type $I^\pm(\gamma)$ are equal for timelike and lightlike curves,
but those type  $\uparrow I^-(\gamma)$, $\downarrow I^+(\gamma)$
are not.}.

(iii) The behaviors of $F$ essential for causality  depend on the
behavior of $F(x,u)$ for large $|x|$ (where $|x|$ denotes the
distance to a fixed point, if $\mo$ is a Riemannian manifold),
\cite{FScqg}. Concretely: (a)  $F$ superquadratic in $x$ and $-F$
subquadratic, which implies that $\m$ is not distinguishing --the
c-boundary is not properly well-defined--, (b) $F$ at most
quadratic, i.e. $F(x,u)\leq R_{1}(u)|x|^{2}+R_{0}(u)$ for
continuous functions $R_i$, which implies that $\m$ is strongly
causal, and (c) $F$ subquadratic and $\mo$
 complete, which implies that $\m$ is globally hyperbolic.

\sm \sm

Notice that there exists a critical behavior  for the causality of
the boundary when $F$ is $x-$quadratic. So, one also defines the
behaviors:

\bit \item $F$ asymptotically quadratic %(in all directions):
$
R^-_{1}(u)|x|^{2}+ R^-_{0}(u) \leq F(x,u)\leq
R_{1}(u)|x|^{2}+R_{0}(u)
$

\item In particular, $F$ is $\lambda$-asymptotically quadratic
($\lambda>0$) when  :
$$\frac{\lambda^{2}|x|^{2}+R^{-}_{0}}{u^{2}+1}\leq F(x,u)\leq
R_{1}(u)|x|^{2}+R_{0}(u) \quad \quad \forall (x,u)\in M\times\R $$
(For the properties of the c-boundary, this behavior will be
enough in some $\mo$-direction.)
 \eit

%\begin{center}
%\includegraphics[width=1.2\textwidth]{Esq1.jpg}
%\end{center}

\subsection{Ingredients for the general computation}

The geometric and analytic tools necessary for the computation are
essentially three:

\sm
\sm

\noindent {\em 1.- Functional approach associated to an arrival
time}. For any $z_0=(x_0,u_0,v_0)\in \m$  parameterize with $u$
each lightlike curve (which is not the integral curve of
 $\partial_v$) starting at $z_0$: $\gamma(u)=(x(u),u,v(u)), u\in I \subset \R$. As
 $\langle \gamma',\gamma'\rangle \equiv 0$ one has:
\be \label{euve}
v(u)=v_{0}+ %\frac{1}{2}
\frac{1}{2}\int_{u_0}^{u}(|\dot{x}(\sigma)|^2-F(x(\sigma),\sigma))d\sigma,
\quad \forall u\in I. \ee
 Consider the lower extreme of $I^+(z_0)
\cap L_{(x_0,u_0)}$, where $L_{(x_0,u_0)}=\{(x_0,u_0,v): v\in
\R\}$. This extreme (minus $v_0$) is the ``arrival time'' from
$z_0$ to the lightlike line $L_{(x_0,u_0)}$. It is equal to the
infimum of all the possible $v(u)$ which can be obtained from
(\ref{euve}).

Now, consider the set ${\cal C} (\equiv {\cal C} (x_0,x_1;|\Delta
u|))$ of curves $x:[0,|\Delta u|]\rightarrow \mo$ joining $x_0$
with some $x_1\in \mo$.
%i.e, of the functional
%$${\cal C}\rightarrow \R, \quad y\mapsto v_y(|\Delta u|).
%$$
%That is,
Then, $I^+(z_0)$ is controlled by the infimum of the ``arrival
time functional'':
$$
 \J_{u_0}^{\Delta u}: {\cal C}\rightarrow \R, \quad
 \J_{u_0}^{\Delta u}(y)=
\frac{1}{2} \int_{0}^{|\Delta
 u|}(|\dot{y}(s)|^2-F(y(s),u_\nu(s)))ds.
$$
where $u_\nu(s)=u_0+\nu (\Delta u) s, \nu =+1$ (and analogously
$I^-(z_0)$).

Notice that $\J_{u_0}^{\Delta u}$ is, in fact, a Lagrangian.  As
we are really interested in, say, $P=I^-(\gamma)$ and $\uparrow
P=\uparrow I^-(\gamma)$, one has to study limits such as:
$$
\hbox{Inf}({\cal J}_{u_0}^{\Delta u}) \; \hbox{ on} \; {\cal
C}(x_0,x_{\Delta};\vert \Delta u\vert) \quad  \hbox{with} \;
x_{\Delta} = x(u_{\Delta}), \; \hbox{and} \; u_{\Delta} =
u_0+\Delta u\nearrow u_{\infty}
$$
Abstract technical conditions on ${\cal J}_{u_0}^{\Delta u} $
(implied by the suitable asymptotic behaviors of $F$ above)
control when the TIPs $P$ ``collapse'' to a 1-dimensional boundary
--as in the examples provided by Berenstein and Nastase, and
subsequent studies.

\sm \sm

\noindent {\em 2.-  Busemann type functions}. In order to check
when pairs $P,F$ are $S$-related, TIP's and TIF's as well as their
common futures/pasts must be computed explicitly.

In standard static spacetimes (conformal to products $\LL^1\times
S$), certain Busemann-type functions, constructed from curves ``in
the spatial part'', were useful to compute TIP's \cite{HaNonlin}.
Here, from curves in $\mo$ one constructs: (a) an adapted version
$b^-$ of such Busemann functions in order to compute  TIP's and
TIF's, and (b) a more refined version $b^+$ in order to compute
common future/pasts.

\sm \sm

\noindent {\em 3.- Sturm-Liouville theory}. In order to study
accurately the limit cases of quadratic Mp-waves (some of them of
special physical interest), the Euler-Lagrange equation for the
functional ${\cal J}_{u_0}^{\Delta u}$ must be analyzed. In fact,
not only this equation must be studied but also the variation of
the solutions with $u_0+\Delta u \nearrow u_\infty$. This yields
an associate Sturm-Liouville problem, especially interesting for
Mp-waves with critical behavior.

\subsection{Table of results}

With all the previous elements, the results can be summarized as
follows. We refer to \cite{FS} for detailed proofs and
explanations.

\begin{center}
\noindent
\begin{tabular}{|c|c|c|c|} \hline
% & & & \\
 \mbox{Qualitative $F$} &  \mbox{Causality} & \mbox{Boundary $\partial M$} & \mbox{Some examples}
  \\
 %& & & \\
 \hline
 %& & & \\
$ \begin{array}{c} \mbox{$F$ superquad.} \\ \mbox{$-F$ at most
quad.}\end{array} $ & $\begin{array}{c} \mbox{No distin-} \\
\mbox{guishing}\end{array}$
 & \mbox{No boundary} & $\begin{array}{c} \mbox{pp-waves yielding} \\
\mbox{ Sine-Gordon string} %\\ \mbox{and related ones}
\end{array}$  \\
 %& & & \\
 \hline
 %& & & \\
 $\begin{array}{c} \mbox{At most quad. $F$} \\
\mbox{(resp.  %$^1$
$|F|$)}\end{array}$
 & $\begin{array}{c} \mbox{Strongly} \\
\mbox{causal}\end{array}$ &  $\begin{array}{c} \mbox{In
principle,}
\\ \mbox{computable}
%\mbox{from Th. \ref{t-pastsets}} \\
%\mbox{(resp. Th. \ref{3})}
\end{array}$ & all below
  \\
% & & & \\
\hline
 %& & & \\
$\begin{array}{c} \mbox{$\lambda$-asymp. quad. %$^2$
} \\
\mbox{$\lambda>1/2$}\end{array}$
 & $\begin{array}{c} \mbox{Strongly} \\
\mbox{causal}\end{array}$ & $\begin{array}{c} \mbox{1-dimension,} \\
\mbox{[lightlike]}\end{array}$ & $\begin{array}{c} \mbox{plane waves} \\
\mbox{with some eigenv.} \\
\mbox{$\mu_1 \geq \lambda^{2}/(1+u^2)$} \\ \mbox{for $|u|$ large }\end{array}$  \\
 %& & & \\
 \hline
 %& & & \\
$\begin{array}{c} \mbox{$\lambda$-asymp. quad.} \\
\mbox{$\lambda \leq 1/2$}\end{array}$
 & $\begin{array}{c} \mbox{Strongly}
\\
\mbox{causal}
\end{array}$ &
$\begin{array}{c} \mbox{Critical}
%\\
%\mbox{dimension}
\end{array}$ &
$\begin{array}{c} \mbox{pp-wave with}
%$\S$\ref{s9.1}}
\\
\mbox{$F(x,u)=\lambda^{2} x^2/(1+u)^2$}\\
\mbox{(for $ u>0$)}
\end{array}$
\\
 %& & & \\
 \hline
 %& & & \\
 Subquadratic & $\begin{array}{c} \mbox{Globally} \\
\mbox{hyperbolic}\end{array}$ & $\begin{array}{c} \mbox{No identif.} \\ \mbox{in $\hat{\partial} \m, \check{\partial} \m$} \\ \mbox{Expected} \\
\mbox{higher dim.}\end{array}$ &
$\begin{array}{c} \mbox{ (1) $\LL^n$ and static} \\
\mbox{type waves}
\\
 \mbox{(2)  plane waves with} \\
\mbox{$-F$ quadratic}
\end{array}$  \\
 %& & & \\
 \hline
\end{tabular} \\
%\vspace{\parskip}

\end{center}

\section{Conclusion}\label{s7}

Finally, let us summarize both, our proposal for the c-boundary,
and how it must be computed ideally --this purges all the
subtleties and approaches reviewed.

The c-boundary makes sense at least for any strongly causal
spacetime $M$ --where all the desired properties for the topology
will hold surely. Now:

\ben \item In order to obtain $\partial M$ as a point set, compute
all the TIP's and TIF's, i.e., $\hat M, \check M$. Then, consider
all the pairs of type: (a) $(P,F)$,  such that $P\sim_S F$
according to S-relation (\ref{eSz}), and (b) $(P,\emptyset),
(\emptyset, F)$ when $P$ or $F$ are not S-related with anything
(for $P\in \hat M, F\in \check M$ always).

Typically, each $P$ or $F$ will be S-related to at most one
element, which will be also rather intuitive (for example, in
standard static spacetimes \cite{AF}, or when  a suitable
conformal boundary exists \cite{FHS}). Nevertheless, wave-type
spacetimes provide examples where the pairing $(P,F)$ is not by
any means evident. And one can construct relatively simple
examples (open subsets of $\LL^n$, or standard stationary
spacetimes \cite{FHSst}), where a TIP or TIF is S-related to more
than one element.

\item Extend the chronological relation to $\bar M$ by using
(\ref{efp}) --of course, regarding any point in $M$ as a pair
$(I^-(p),I^+(p))$.

The timelike part $\partial_0M$ of the boundary  $\partial M$ is
composed of the pairs $(P,F)$ with $P\neq \emptyset \neq F$. It is
empty if and only if $M$ is globally hyperbolic. The points with
$F=\emptyset$ (resp. $P=\emptyset$) define the future infinite
$\partial_+M$ (resp. the past infinity $\partial_-M$) of $M$.

\item The topology in $\bar M$ is defined according to Definition
\ref{dtop} (see also (\ref{etop}), (\ref{dec}), (\ref{limitop})).
 For a sequence
$\{Q_n=(P_n,F_n)\}_n$ in $\bar M$ with $LI\{P_n\}=LS\{P_n\}=:
P_\infty$, $LI\{F_n\}=LS\{F_n\}=: F_\infty$, this will mean the
following. A point $Q=(P,F)\in \bar M$ lies in the closure of
$\{Q_n\}_n$ iff either $P=P_\infty$ and $F=F_\infty$ or, at least,
$P, F$ are maximal among all the (past or future) indecomposable
sets included in $P_\infty, F_\infty$, respectively --that is:
$P\subset P_\infty, F\subset F_\infty$ and no $P'\in \hat M, F'\in
\check M$ satisfy either $P \subsetneq P' \subset P_\infty$ or
$F\subsetneq F' \subset F_\infty$.%\footnote{Confirmar y discutir.
%P. ej., si basta con considerar sucesiones como las de arriba para
%determinar la topologia -lo cual seria comodo...}.

Then, all  the reasonable properties of a boundary hold for
$\partial M$ (it is closed in $M$, points of the boundary are
$T_2$-separated of the points of the spacetime $M$, etc.). In
particular, all the points in $\partial M$ are $T_1$ separated,
and the cases where they are not $T_2$ separated are intuitively
clear and acceptable.

 \een

{\small

 \edoc